\documentclass[12pt,a4paper]{article}
\usepackage{amsmath, amsthm, amssymb, mathrsfs, amscd, amsthm, amscd,amsfonts,graphicx}
\usepackage[T1]{fontenc}
\usepackage[utf8]{inputenc}
\usepackage{epsfig}
\usepackage{faktor}
\usepackage[all]{xy}
\usepackage[german]{babel}
\newcommand{\ko}{{\mathcal O}}

\newcommand{\kh}{{\mathcal H}}
\newcommand{\kf}{{\mathcal F}}
\newcommand{\Z}{\mathbb{Z}}

\renewcommand{\P}{\mathbb{P}}

\newcommand{\C}{{\mathbb C}}
\newcommand{\R}{{\mathbb R}}
\newcommand{\SL}{{\text{SL}}}
\newcommand{\GL}{{\text{GL}}}

\begin{document}
\vspace{1cm}

\begin{center}
\textbf{\Large Life and work of Egbert Brieskorn (1936 -- 2013)\footnote{Translation of the German article ''Leben und Werk von Egbert Brieskorn (1936 -- 2013)'', Jahresber. Dtsch. Math.--Ver. 118, No. 3, 143-178 (2016).}}


\textbf{\em  Gert-Martin Greuel, Walter Purkert}
\end{center}
\bigskip
\bigskip

\begin{center}
\includegraphics[width=7cm,height=8cm]{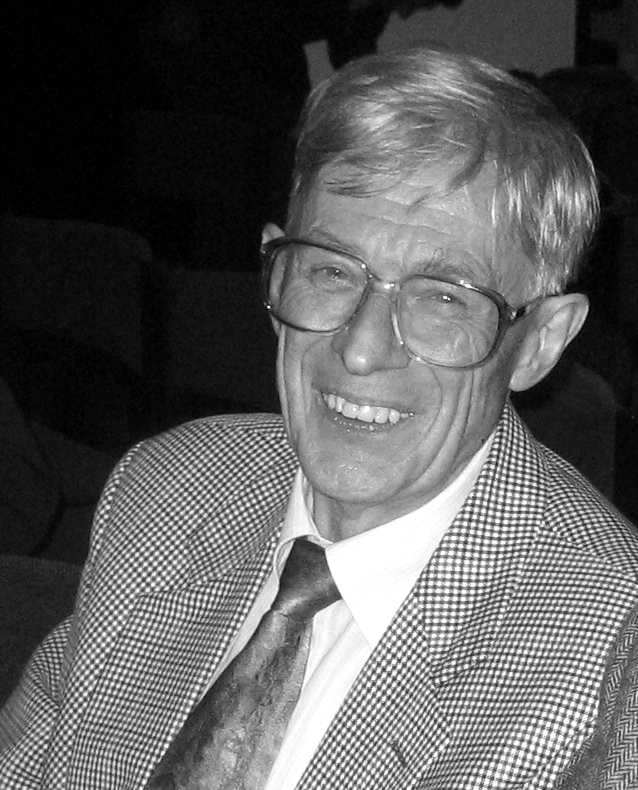}\\
Brieskorn 2007
\end{center}

\bigskip

Egbert Brieskorn died on July 11, 2013, a few days after his 77th birthday. He was an impressive personality who left a lasting impression on anyone who knew him, be it in or out of mathematics. Brieskorn was a great mathematician, but his interests, knowledge, and activities went far beyond mathematics. In the following article, which is strongly influenced by the authors' many years of personal ties with Brieskorn, we try to give a deeper insight into the life and work of Brieskorn. In doing so, we highlight both his personal commitment to peace and the environment as well as his long--standing exploration of the life and work of Felix Hausdorff and the publication of Hausdorff 's Collected Works. The focus of the article, however, is on the presentation of his remarkable and influential mathematical work.

The first author (GMG) has spent significant parts of his scientific career as a graduate and doctoral student with Brieskorn in G\"ottingen and later as his assistant in Bonn. He describes in the first two parts, partly from the memory of personal cooperation, aspects of Brieskorn's life and of his political and social commitment. In addition, in the section on Brieskorn's mathematical work, he explains in detail the main scientific results of his publications. The second author (WP) worked together with Brieskorn for many years, mainly in connection with the Hausdorff project; the corresponding section on the Hausdorff project was written by him.

We thank Wolfgang and Bettina Ebeling, Helmut Hamm, Thomas Peternell, Anna Pratoussevitch and Wolfgang Soergel for useful information and especially Brieskorn's wife Heidrun Brieskorn for the release of material from Brieskorn's estate.
We also thank Andrew Ranicki for encouraging us to translate the article into English and 
special thanks to him and Ida Thompson for checking the translation.
\bigskip

\textbf{Mathematical Subject Classification (MSC2010)}:  01A61, 14B05, 14B07, 14D05, 14F45, 14H20, 14J17, 14J70, 17B22, 32S05, 32S25, 32S40, 32S55, 57R55\\

\textbf{\large Stations of his life}

Brieskorn was born on July 7, 1936 in Rostock, Germany, the son of a mill--construction engineer, and grew up with his sister and his mother in Siegerland. Little is known about his youth and the source of his enthusiasm for mathematics. But from the chapter {\em Childhood and Education}
from the Simons Foundation film about Brieskorn \cite {EB2010} we know that his mother supported his childlike curiosity and that his father promoted his technical interest. He also had a good maths teacher in grammar school, who provided him with mathematical literature beyond the subject matter of the
curriculum. He was particularly interested in geometric constructions (and less, for example, in a work by Gauss).

Even though his technical interest initially prevailed, his interest in mathematics was already strong before his studies. At the examination for admission to the Evangelische Studienwerk Villigst, the funding organization of the Protestant Church for gifted students, the examiner said to him: `Mr Brieskorn, your talent and your enthusiasm for mathematics are extraordinary but do not forget that there are other things besides mathematics in life.'   Egbert Brieskorn told this episode to the first--named author of these lines with a slightly ironic undertone much later, when in fact other things than mathematics  determined his life and work. The head of the Evangelische Studienwerk recognized that his original desire to study electrical engineering was not right for his nature and convinced him to study something theoretical.

Brieskorn therefore began to study mathematics and physics in Munich in October 1956. After five semesters he followed the advice of Karl Stein and moved to Bonn for the summer semester of 1959, in order to understand the theorem of Hirzebruch-Riemann-Roch, which he described as `my first love in mathematics' \cite{EB2000}. Friedrich Hirzebruch, who himself had only come to Bonn in 1956, deeply impressed the young student Brieskorn with his friendly, open personality and his clear style of presentation. Brieskorn became Hirzebruch's student and received his doctorate in 1963 with the thesis ''Differentialtopologische und analytische Klassifizierung gewisser algebraischer Mannigfaltigkeiten'' (Differential topological and analytical classification of certain algebraic manifolds). Hirzebruch later described Brieskorn as his most talented student and Brieskorn highly revered his teacher Hirzebruch as a mathematician and human being all his life.

In 1968, Brieskorn habilitated in Bonn with the thesis {\em Singularit\"aten komplexer R\"aume} (Singularities of complex spaces) and was appointed full professor in G\"ottingen in 1969, where he remained until 1973. Because of his wife Heidrun, whom he married in 1973 and who got a job as a violist at the Cologne Radio Symphony Orchestra (today the WDR Symphony Orchestra Cologne), he moved in 1973 to Bonn, first to the Sonderforschungsbereich Theoretische Mathematik and from 1975 to a position as a full professor, where he worked until his retirement in 2001.

\begin{center}
\includegraphics[width=12cm,height=7.5cm]{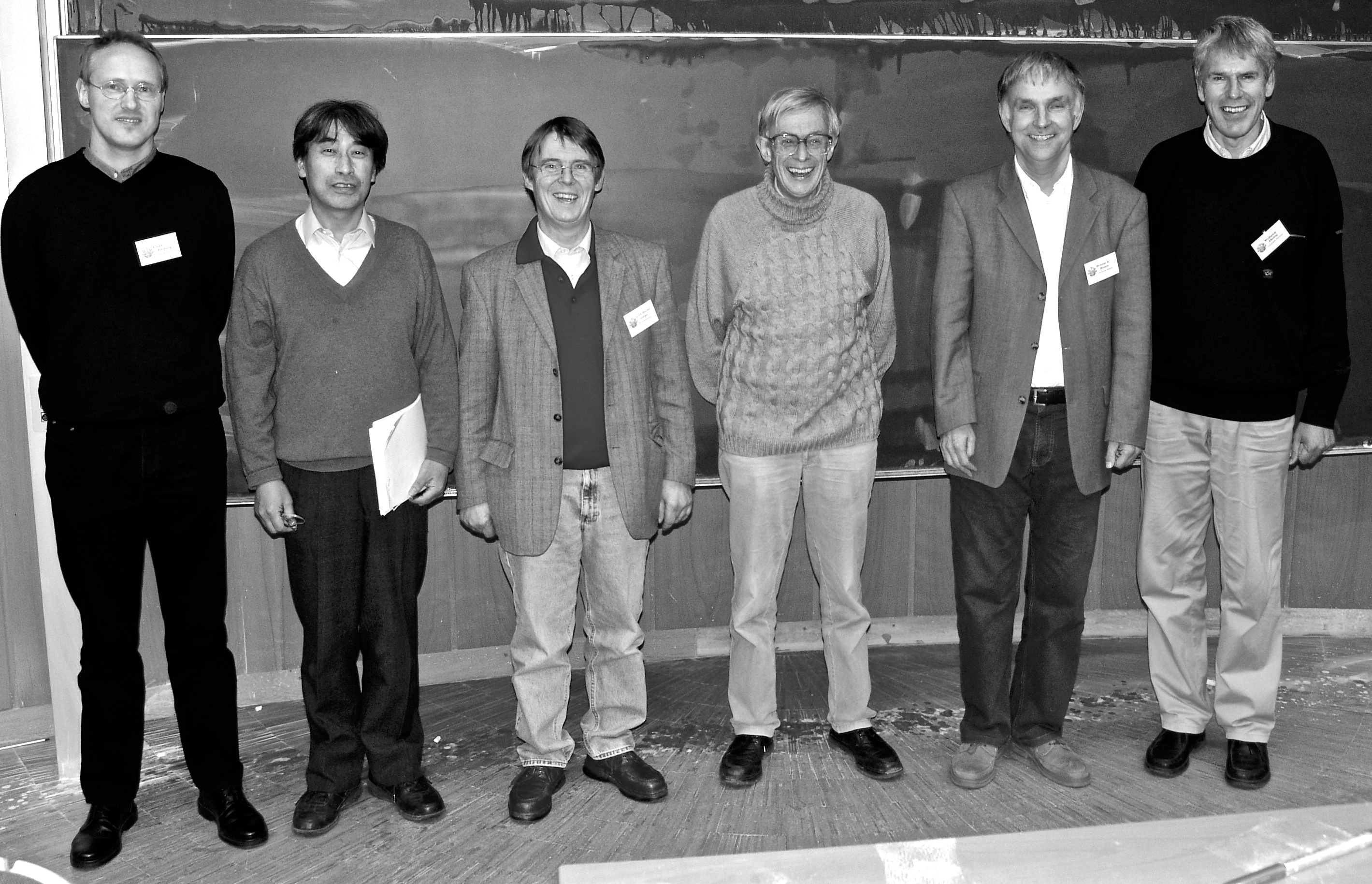}\\

Brieskorn (3. from r.) and his students (from l.) Claus Hertling, Kyoji Saito, Gert-Martin Greuel, Helmut Hamm, Wolfgang Ebeling, 2004
\end{center}

\bigskip

\textbf{\large Mathematics and political--social engagement}

Although Brieskorn has always been socially engaged, e.g. by working in a steel rolling mill of the ``Dortmunder H\"order H\"uttenunion'' in the context of a work semester of the Evangelischen Studienwerk, for him, during his studies and many years as a professor, mathematics was the most important thing in his life, as he writes himself (curriculum vitae work semester). He was completely excited by the fascinating beauty and clarity and the high standards of mathematics.
But he saw mathematical phenomena everywhere, also in small things, and was fascinated by them.
This enthusiasm for mathematics and his own field of research, the theory of singularities, was passed on to the students.

The following sections show the development of the relation of mathematics and political--social engagement of Brieskorn from the perspective and memory of the first-named author.
In the late summer of 1969 he had invited me to tea in his apartment in G\"ottingen to discuss a topic for my diploma thesis. Before he started, he discovered  a refractory caustic in his teacup, which he classified as a ``simple singularity'', and then interpreted the spiral chocolate trail in the biscuits as a dynamical system. He enthusiastically explained to me exotic spheres and how they can be described by real-analytic equations as the boundaries of certain isolated hypersurface singularities.
I was infected, and when he then suggested that I generalize the results of his unpublished work, {\em Die Monodromie der isolierten Singularit\"aten von Hyperfl\"achen} (The monodromy of isolated singularities of hypersurfaces) \cite{EB1970a}, to isolated singularities of complete intersections, I immediately accepted.

I became a student of Brieskorn on the recommendation of Hans Grauert, whose beginner course ``Differential and integral calculus'' I had attended in G\"ottingen in 1966/67. After returning from a one--year study visit at ETH Zurich, I received a telephone call from Brieskorn asking if I would be interested in writing a thesis with him. Grauert recommended me because we were in the same ``fraternity'' (i.e. the Evangelische Studienwerk). So it happened that I became Brieskorn's first graduate student.

Brieskorn had come to G\"ottingen in the summer semester 1969, first as
a substitute professor and from July 1969 as a full professor. In the winter semester 1969/70 he was on leave for a research stay at the IHES in Bures-sur-Yvette. I had missed his first lecture in the summer of 1969 in G\"ottingen on 2-dimensional schemes and over winter I learned sheaf theory and hypercohomology in the reading room of the Mathematical Institute in the Bunsenstrasse with the help of Godement's {\em Topologie
Alg\'{e}brique et Th\'{e}ory des Faisceaux}. Except for occasional whispering and the unmistakable sound of the heavy breathing of Carl Ludwig Siegel, who  disappeared into the back of the reading room to consult the older works, one was undisturbed and could discuss with fellow students in the adjoining discussion room.

After returning from France, Brieskorn gave the following lectures: ``Differential Topology'', then the beginner lecture ``Calculus I and II'', then ``Analysis on Manifolds'', ``Qualitative Theory of Dynamical Systems'', ``Algebraic Topology II'', and before he moved to Bonn in the summer semester 1973 ``Simple Singularities''. Brieskorn's lectures differed fundamentally from those of Grauert. While Grauert only told what he wrote and proved on the blackboard, Brieskorn often presented larger connections and mentioned cross-connections without proving them. Obviously, both approaches have their advantages, and the students in G\"ottingen greatly appreciated both Brieskorn's and Grauert's lectures.

Brieskorn often made an extraordinary effort in the preparation of his lectures, to show mathematical and historical backgrounds or side branches of the  material he treated. This can be clearly seen in his textbooks. It can be said that his quest for perfection, which even increased over the years, was characteristic of him. He also expected perfection from his students, which he greatly promoted, both mathematically, among others with weekly working meetings or with recommendations, but also in the personal sphere. For example, he let his first doctoral student Helmut Hamm live in his apartment during his stay in France free of charge.

The G\"ottingen period from 1969 to 1973 caused a remarkable change in Brieskorn's views and attitudes. In addition to mathematics, political and social issues gained importance for him. It was the time shortly after the violent student protests of the '68 generation, in which the students of G\"ottingen, primarily the theologians, followed by the mathematicians, were very active (for example in a left-wing action group ``Basisgruppe Mathematik''). Brieskorn was very positive about some of the demands of the students and rather critical of others. He objected to scientifically unfounded hierarchical structures, was committed to greater co-determination of research assistants and students, and he sympathized with ideas of the reform universities in Bremen and Osnabr\"uck. However, the scientific quality of the study always came first.

Even more important than student reform ideas for Brieskorn, however, were his commitment against the Vietnam War and, in turn, his commitment to supporting oppressed peoples. From [26] we know that during his stay at the Massachusetts Institute of Technology (MIT) he participated in a major demonstration against the Vietnam War, together with Michael Artin and other MIT colleagues in New York, at which Martin Luther King spoke.

From the beginning he was involved in G\"ottingen  in the Committee for Scientific Cooperation with Cuba (KoWiZuKu), which was founded in 1970 and whose first Secretary General was the mathematician Klaus Krickeberg. The fact that he was always very precise in these activities is shown by an episode when, together with me and my wife in G\"ottingen, he stuck up posters announcing a lecture by the Cuban health minister at the university. Brieskorn meticulously ensured that everything was correct, and no posters were glued to distribution boxes as this could lead to heat problems.

Later in Bonn Brieskorn was involved from the outset in the peace movement, which was formed in protest against the 1979 ``NATO double--track decision'' and thus against the deployment of medium-range missiles in the then FRG. Brieskorn was one of the initial signatories of the ``Mainzer Appell'', the final declaration of the congress ``Responsibility for peace -- scientists warn against new nuclear armaments'' in June 1983 in Mainz, in which the mathematician Stephen Smale also participated (see \cite {HD1983}). Especially physicists, but also many mathematicians such as Brieskorn, organized themselves in the ``Scientist Initiative Responsibility for Peace'' (today: ``Scientist Initiative for Peace and Sustainability'').

Brieskorn's political views in the Bonn era can certainly be classified as left of social-democratic ideas, but over the years they have evolved into radical ecological convictions, which he supported together with his wife Heidrun. They have lived since 1982 in a house at an isolated location on the edge of the forest in Eitorf an der Sieg where both, together with a small group, have devoted themselves intensively to nature conservation and more specifically to the conservation of species. The species in question were initially indigenous bats, for which they controlled winter quarters in old mine tunnels and secured and built new winter quarters. In order to distinguish the different species, Heidrun Brieskorn made many sound recordings, which Egbert Brieskorn subjected to self-written programs using a Fourier analysis.

Both devoted even more time and effort to protect and preserve the living conditions of a very rare species of butterfly, the large blue (Latin: {\it Maculinea}). This went so far that Brieskorn and his wife persuaded the community to change their development plans and they themselves bought up grassland to maintain the habitat of these butterflies. They founded the Maculinea Foundation NRW so that the work to preserve the butterfly species can continue on a permanent basis. For their commitment they were jointly awarded in 2013 the decoration ``Bundesverdienstkreuz am Bande''  (``Cross of the Order of Merit'') of the Federal Republic of Germany.

In addition to the volunteer work in nature conservation, the last 20 years of Brieskorn's life were determined by his collaboration in the edition project ``Felix Hausdorff -- Collected Works'' of the North Rhine-Westphalian Academy of Sciences and the Arts. Brieskorn himself wrote shortly before his 75th birthday in June 2011, knowing that he might not live much longer, in a letter to `my dear former students and my students' how it came about and how much he was concerned with the biography of Hausdorff:
``One of the tasks developed from the fact that the Mathematical Institute in Bonn in January 1992 wanted to celebrate the 50th anniversary of the death of Felix Hausdorff. Since no colleague wanted to give a lecture about his life, I took over this task at that time, not knowing what I had let myself in for. I have spent 20 years searching for archives and sources of all kinds for traces of this extremely remarkable man and mathematician. I learned a lot while sacrificing a lot of time. I made three attempts to write his biography, and with the third version, I think, I am on the right track. About 262 pages are written so far, but not even half of the path of his life is described. This biography is to appear in the first volume of a total of nine-volume edition of Felix Hausdorff's works. Six volumes are printed, and at least one or perhaps two of the missing volumes will be printed this year. The bad thing is that the first volume in which my biography of Hausdorff is to be released cannot possibly be completed on time.''

The following section reports on Brieskorn's extensive work and his research on this project, which goes far beyond the usual.
\bigskip

\textbf{\large The Hausdorff project}

The work, which Egbert Brieskorn mentions in his letter, was a labour of love for him for more than 20 years, in particular the research into the life and work of Felix Hausdorff (1868-1942) and the publishing of Hausdorff's collected works.

Felix Hausdorff founded general topology as a freestanding area of mathematics and,  in addition, made fundamental contributions to general and descriptive set theory, measure theory and analysis. His contributions to the theory of Lie algebras, probability theory and actuarial mathematics were also important to subsequent developments. Hausdorff was also (for a mathematician somewhat uniquely) a notable writer and philosophical author. From 1897 until 1913, under the pseudonym of Paul Mongr\'e, he published a volume of aphorisms, an epistemological book, a volume of poetry, a play which has been performed more than 300 times in over 30 cities, as well as 17 essays in leading literary journals. In the twenties and early thirties, he was internationally recognised and respected as the leading representative of the Mathematical Institute of Bonn. As a Jew living under the dictatorship of the National Socialist Party, he and his wife took their own lives on 26th January 1942, when deportation to a concentration camp was imminent.

   Egbert Brieskorn's involvement began in 1979, when the suggestion came from the students of the Mathematical Institute to honour Felix Hausdorff through a memorial plaque in the institute. Brieskorn, who had already been involved in the peace movement for years and who had talked in depth about questions concerning the social responsibility of scientists, supported the plan from the very beginning, gave his own financial support and organised a collection of donations amongst the teaching staff. In 1980, on the anniversary of Hausdorff's death on 26th January, the marble plaque was unveiled at the old institute on Wegelerstrasse. It was there for over 30 years, until it was recently
   transferred to the new institute building. On the occasion of the unveiling, an article appeared about Hausdorff's tragic fate penned by the historian Herbert Mehrtens. Egbert Brieskorn wrote an introduction for the article, the closing sentences of which are quoted here; he wrote: ``No form of inhumanity and oppression must leave us indifferent, even when the victims are distant und unfamiliar to us. The thought of a person like Hausdorff, whom we all recognise for his great scientific achievements, can also help us to sharpen our conscience and sense of responsibility. Without a growing feeling of responsibility, scientists will not be in the position to make their contribution to a more humane society. For this task we also need to come to terms with the past.''

In November 1980, Egbert Brieskorn was successful in procuring Hausdorff's literary estate for the university library in Bonn, consisting of some 26,000 pages in size and since 1964 in the private ownership of Prof. G\"unther Bergmann in M\"unster. The proceeds benefited Hausdorff's daughter Lenore K\"onig, who lived in poor conditions in Bonn in a retirement home.  As the contract was signed, Egbert Brieskorn wrote to G\"unther Bergmann on 15th November: ``The arrangements, which are now being made, seem to me to be very good. It is true that I am no historian, but I do have a particular interest in the history of mathematics and occasionally oversee dissertations with historical aspects. I hope that one day there will be a mathematical historian who will work on Hausdorff. Then the literary estate in our university library will be very important.'' Whilst writing these lines, he had surely not thought that it would be he himself who would 10 years later tackle the project of writing Hausdorff's biography.

26th January 1992 was the 50th anniversary of Hausdorff's death. On this occasion, Egbert Brieskorn organised a memorial colloquium of the Mathematical Institute of the University of Bonn, the result of which was the volume {\em Felix Hausdorff zum Ged\"achtnis  -- Aspekte seines Werkes} (In Memory of Felix Hausdorff -- Aspects of His Work), edited by Brieskorn himself and published by Vieweg-Verlag. Furthermore, he arranged an exhibition about the life and work of Felix Hausdorff, which was met with a lively interest not just amongst the Mathematical Institute, but also the university and the general public of Bonn. For the exhibition, he published a comprehensive catalogue with an initial biographical sketch of Hausdorff. On 1st February 1992, there was an hour--long radio programme for the series ``Mosaik'' (Mosaic) on WDR, about Hausdorff and the exhibition in Bonn under the title ``Auf d\"unner Schneide tanzt mein Gl\"uck'' (My Happiness Dances on the Edge of a Knife), which is the first line of Hausdorff's poem ``Wiederkunft'' (Return) from the poetry volume ``Ekstasen'' (Ecstasy). Egbert Brieskorn and the WDR culture editor Friedrich Riehl conceived the programme together, and made the exhibition and its subject well known far beyond Bonn. Even the media coverage of the exhibition was significant.

The preparation for all of these activities began in 1989 with numerous conversations, which Egbert Brieskorn held with Hausdorff's daughter, with contact with Hausdorff's niece Else Pappenheim and with further contemporary witnesses, as well as with the collection of material for the biography, in particular through researching in archives. During this research, he showed thoroughness and the mind of a detective, which any professional historian would hold in high esteem. In his literary estate there are dozens of thick folders, documenting all of these endeavours and in successful instances, though that was not always the case, recording their results.

In this spirit and in the run-up to the memorial colloquium, the mathematicians from Bonn, Egbert Brieskorn, Friedrich Hirzebruch and Stefan Hildebrandt, discussed the possibilities and necessary steps for setting in motion an edition of the works of Hausdorff, under the consultation of external experts. Friedrich Hirzebruch suggested creating a Hausdorff Commission at the North Rhine-Westphalian Academy of Science, which was already overseeing a number of edition projects, to plan and then supervise a similar project. The academy agreed and the commission took up its work under the leadership of academy member Reinhold Remmert in 1991. As the first step, a careful indexing and cataloguing of the Hausdorff literary estate was envisaged. This  was carried out from October 1993 until the end of 1995. The result was an inventory of 550 pages with short descriptions of the content of each individual fascicle. After this vital step towards the success of the edition project was completed, and with the inclusion of the literary estate, the final project was within sight.

In order to create a diligently commented edition with the inclusion of selected parts from the literary estate in the work, many things needed to be considered and done.  One must establish editorial principles, develop a volume structure, find suitable personnel and, crucially, applications must be made in order to finance the project. Egbert Brieskorn was the guiding spirit throughout all of these discussions and activities. Particularly difficult was the winning over of suitable philosophers and literary academics as editors of the volumes dedicated to that side of Hausdorff's creations. For this purpose, he took part in a philosophical conference about Jewish Nietzscheanism in Greifswald, in order to come into contact with appropriate experts, and in Bonn organised an interdisciplinary conference for academics in the humanities, with the support of the MPI for Mathematics, about Hausdorff's philosophical and literary work.

The application for the financing of the project, which he submitted to the DFG in 1996 together with Friedrich Hirzebruch and Erhard Scholz, was finally approved and in November of 1996 the working team ``the Hausdorff--Edition'' at the Mathematical Institute took up work under his leadership. In 2002, the North Rhine-Westphalian Academy of Sciences and Arts took over the Hausdorff--Edition as one of their academy--projects. Some of the originally employed editors were not able to work on the project due to various reasons, meaning that new employees frequently needed to be found.  In the end 14 mathematicians, four mathematical historians, two literary academics, one philosopher and one astronomer from four countries collaborated on the editing.

Egbert Brieskorn was particularly concerned during the editing with tracking the four noticeable threads, not immediately obvious, which lead from Hausdorff the mathematician to Mongr\'e the writer and philosopher.  For this reason he always emphasized the inter--disciplinary character of the project, holding pioneering lectures at three large editorial conferences exactly in this respect and organising a series of discussion circles, which involved mathematicians, philosophers and literary academics. In order to obtain an impression of his intentions, we quote here the beginning of his programmatic lecture at one of these conferences in February 2003 at Schloss Rauischholzhausen:  ``\,\grq It is not always determined that the concept must lie within the lines or indeed only between the lines, perhaps it is somewhere else entirely, far, far away! Perhaps the author sounds his bell, and somewhere a string with the same vibration and the same tone colour answers -- and it is not the actual bell, but the sound of the string that expresses the original concept.\grq  \ We should always keep in mind the sense of this aphorism from Paul Mongr\'e during our round-table discussion about this author, about the author Mongr\'e and about the mathematician Felix Hausdorff. Each one of us will hear something different from the variety of motifs and the richness of tone colour. What someone hears is highly dependent on the hearing experience, which one has made in the course of his or her life. If literary academics, philosophers, historians and mathematicians also listen to each other, we can hope here and there to hear an original concept.''

As leader of the working group, Egbert Brieskorn allowed a wide latitude to the co--ordinator of the edition project and the editors and employees of the individual volumes. When problems arose, he helped with advice. He did however reject some of the submissions when they did not satisfy his high demands. In these cases, decisive improvements could always be made. In the meantime, eight of the ten planned volumes were released with Springer: volume IA {\em General Set Theory} (2013), volume II {\em Basics of Set Theory} (2002), volume III {\em Descriptive Set Theory and Topology} (2008), volume IV {\em Analysis, Algebra and Number Theory} (2001), volume V {\em Astronomy, Optics and Probability Theory} (2006), volume VII {\em Philosophical Works} (2004), volume VIII {\em Literary Works} (2010), and volume IX {\em Correspondence} (2012). Egbert Brieskorn's voice is perceptible in all of the volumes, even when he did not work explicitly on every volume himself.

 He did, however, take on the most difficult part of the project himself: volume IB, the biography. Here, next to the mathematical work of Hausdorff, were also further fields of interest and references to his life from very different areas: philosophy, in particular Kant, Schopenhauer, Nietzsche and Hausdorff's relationship to the Nietzsche archive, epistemology, in particular Hausdorff's language criticism and his thoughts on literary figures like Dehmel, Hartleben and Wedekind, music, particularly Hausdorff's relationship to Wagner and his relationship to Reger, and graphic art, particularly Hausdorff's friendship with Max Klinger. In this volume were also the family history in terms of the Jewish story and the history of anti--Semitism, up until Hausdorff's tragic end under the Nazi dictatorship.

In 2007 Egbert Brieskorn, Erhard Scholz and the co-author of this obituary had the opportunity to present the Hausdorff edition project in the S\'eminaire d'Histoire des Math\'ematiques de l'Institut Henri Poincar\'e in Paris. There Brieskorn gave the introductory speech, in which he said the following about his own work retrospectively (he spoke French of course): `With regard to my own portion of the project, I must confess now that I am no historian and that I am not led predominately by a historical interest. On the contrary, my personal interest sprang originally from two motives. One of the motives was shame about the terrible guilt that Germany, through the persecution and murder of Jews in Europe, has brought upon itself. The other motive was very personal: In the eighties I met Felix Hausdorff's daughter Lenore K\"onig, who at that time was living in a retirement home in Bonn. She gave me an initial introduction to the life and personality of her father. Through this I later felt the personal obligation to better understand the life of this unique person. When the university of Bonn prepared a memorial event for Hausdorff's fiftieth anniversary of death, and later as the plan took shape for an edition of his works, I saw the biography of Felix Hausdorff as my personal task. At the beginning, I underestimated this task concerning the difficulties and also in respect to what this task meant to me personally. This work has greatly changed my own life and way of thinking: I have -- or at least I hope -- learnt something from Hausdorff.'

Egbert Brieskorn worked intensively on the transcript of the biography in his last years, also during his severe illness. Three weeks before his death he sent the last sub--chapter that he had still been able to complete to the working group in Bonn, it was a particularly difficult chapter about the relationship between Hausdorff and the philosopher, mystic and anarchist Gustav Landauer. To this day there are 546 pages of the biography written by Brieskorn and ready for printing. When he sensed that he could not manage any more, he proposed a meeting, in order to explain how he had imagined the rest of the biography. He suggested 12th July 2013 as the date for the meeting. On the evening of 11th July he passed away. It is a duty for the working group in Bonn to bring the volume to completion as well as we can. There are over 100 folders with material for the biography, which he had collected over more than 20 years through often painstakingly detailed work. They were given to us by
Ms Brieskorn, and are of invaluable help.
Volume IB is now finished and will appear at Springer in spring 2018, in the year of Hausdorff's
150th birthday.
\bigskip
\medskip

\textbf{\large Brieskorn's mathematical work}

The mathematical work of Brieskorn is largely determined by the development of the ``singularity theory'' of complex hypersurfaces. In particular, his early work had a great influence on the development of singularity theory and it is no exaggeration to call Brieskorn, along with Vladimir Igorevich Arnold, John W. Milnor and Ren\'e Thom, one of the fathers of singularity theory. In his review lecture \cite {EB1976} Brieskorn writes:

``Singularities exist in all areas of mathematics and in many applications, and the pair of opposites 'regular -- singular' is one of the most common in a whole series of such opposite pairs in mathematics. What is meant by singularities is shown by the analysis of the many different definitions of singular objects. Such an analysis leads to a few basic meanings: a singularity within an entity is a place of uniqueness, peculiarity, degeneration, indefiniteness or infinity. All these meanings are closely related.''

I use the term ``singularity theory'' here in the sense of exploring systems of finitely many differentiable, analytic or algebraic functions near a point in which the Jacobian matrix of the functions does not have maximum rank. This implies, according to the implicit function theorem, that the zero set of the functions is not a differentiable, analytic or algebraic manifold at a singular point. Here, singularities of vector fields or differential forms are included.

The term singularity theory was, to my knowledge, introduced by V.I. Arnold, although it is not really a closed theory with more or less uniform methods.
On the contrary, a characteristic of this field is the variety of different methods used, and the relationships to many other mathematical disciplines. These include algebraic geometry, complex analysis, commutative algebra, combinatorics, representational theory, Lie groups, topology, differential topology, dynamical systems, symplectic geometry, and others. Brieskorn was  particularly fascinated by the complexity and the manifold interactions of singularity theory with other mathematical and non-mathematical fields, such as  geometric optics, and, as we shall see, he has contributed significantly to the study of some of these interactions. However, he has never been able to make friends with the term singularity ``theory''.

Almost all of Brieskorn's mathematical works either deal directly with singularities of complex hypersurfaces, or they are motivated by the study of these singularities. His work shows, in addition to originality and depth, a wide range of questions and methods, which are typical of the entire area.

In the following review of Brieskorn's work, I also try to highlight important results on which Brieskorn's works are based, as well as the developments that resulted from his work. A short description of the scientific work of Brieskorn can be found in \cite {GG1997}.
\bigskip

\textbf{Dissertation}

Singularities do not play any role in the first two publications of Brieskorn that are parts of his dissertation \cite {EB1962}, which he wrote in 1962 as a student of Friedrich Hirzebruch and which he published in \cite {EB1964} and \cite {EB1965}.

The question there is to what extent does the differentiable structure of a K\"ahler manifold already determines its biholomorphic structure in the case of the complex quadric $ Q_n $ or the holomorphic $ \P^n $ bundles over $ \P^1 $. This problem had been studied and solved by F. Hirzebruch and K. Kodaira in 1957 for the complex projective space $ \P^n $.

The following main result of the first paper is an exact analogue of the mentioned theorem of Hirzebruch and Kodaira. Brieskorn had already announced the result in 1961 in the Notices of the AMS:

\textit {Let $ X $ be a $ n $-dimensional K\"ahler manifold that is diffeomorphic to the $ n $-dimensional complex projective quadric $ Q_n $. Then:}
\textit{
\begin{enumerate}
\item [(i)] If $ n $ is odd, then $ X $ is biholomorphic to $ Q_n. $
\item [(ii)] If $ n $ is even and $ n \neq  2$, then the 1st Chern class $ c_1 $ of $ X $ satisfies: $ c_1 = \pm ng $, where $ g $ is the positive generator of $ H ^ 2 (X, \Z) \cong \Z$; if $ c_1 = + ng $, then $ X $ is biholomorphic to $ Q_n $.
\end{enumerate}}

Brieskorn asks if the assumption that $ X $ is K\"ahler can be omitted and if there are any K\"ahler manifolds with $ c_1 = -ng $ that are diffeomorphic to $ Q_n $. Both problems seem to be open to this day.

As a corollary Brieskorn proves statements about the deformation behavior of $ Q_n $. He considers a family of complex manifolds $ V_t, \ t \in \C, \ | t | $ sufficiently small, and proves:
\textit{
\begin{enumerate}
\item [(i)] If $ V_0 \cong Q_n $, then $ V_t \cong Q_n $ \ for $ t \neq 0 $,
\item [(ii)] If $ V_0 $ is  K\"ahler and $ V_t \cong Q_n$ \ for $ t \neq 0 $ and $ n \geq 3 $, then $ V_0 \cong Q_n $.
\end{enumerate}}

Again, he asks if in (ii) the assumption that $ V_0 $ K\"ahler can be omitted.

That this is indeed the case was proved by J.M. Hwang in 1995, after the same question of \ ``non-deformability'' of $ \P ^ n $ instead of $ Q_n $,  had previously been answered positively by Y.-T. Siu.

For $ n = 2 $ the first statement does not hold, because on the differential manifold $ Q_2 $  there are, according to Hirzebruch, infinitely many different complex structures, the so-called Hirzebruch $ \Sigma $--surfaces $ \Sigma_ {2m} $.

The $ \Sigma $--surfaces are total spaces of holomorphic fiber bundles over $ \P ^ 1 $ with fiber $ \P ^ 1 $. In the second part of his dissertation, Brieskorn examines holomorphic fiber bundles over $ \P ^ 1 $ with fiber $ \P ^ n $, which he calls \textit {$ \Sigma $--manifolds} in accordance with the Hirzebruch $ \Sigma $--surfaces. Taking advantage of Grothendieck's splitting theorem for vector bundles over $ \P ^ 1 $, Brieskorn classifies the $ \Sigma $--manifolds up to biholomorphic and birational equivalence and up to diffeomorphism. As a result, as in the case of $ \Sigma $--surfaces, there is an infinite number of different complex structures on every differentiable $ \Sigma $--manifold. In addition, he proves that $ \Sigma$--manifolds deform into $ \Sigma$--manifolds and that in a K\"ahler family of $ \Sigma $--surfaces, they specialize in a $ \Sigma $--surface (similar to (ii) above).

The questions and the methods of proof in his dissertation come from algebraic and analytic geometry, as well as from algebraic topology. These methods, including sheaf theory, which came from France, were quite new at that time and began to slowly gain acceptance in Germany, especially in the generation of young mathematicians. In an exemplary way they were embodied by Brieskorn's teacher Friedrich Hirzebruch. In addition to Hirzebruch, Hans Grauert and Reinhold Remmert also had great influence on the development of modern analytic and algebraic geometry in Germany, in particular  with the development of the theory of general complex spaces, whose structure might contain nilpotent elements. Brieskorn also thanked Reinhold Remmert and Antonius van de Ven in his dissertation.  He met both of them while working as a research assistant and employee in Erlangen in 1962, showing that he had an extremely inspiring environment for modern mathematics. He was decisively influenced by the spirit of optimism that prevailed in Germany at that time, and especially by the support of his teacher Hirzebruch.
\bigskip

\textbf{Deformation theory}

Friedrich Hirzebruch's book {\em Neue topologische Methoden in der algebraischen Geometrie}, published as early as 1956 in the Springer series ``Ergebnisse der Mathematik und ihrer Grenzgebiete'', had just been published in the second, extended edition. The theorem of Hirzebruch--Riemann--Roch, which was proved there, was one of the foundations of Brieskorn's dissertation. Another basis was the deformation theory of analytic structures developed by K. Kodaira and D. C. Spencer.

The Hirzebruch--Riemann--Roch theorem was a great generalization of Riemann--Roch's classical theorem to complex vector bundles on arbitrary complex projective manifolds, rather than divisors on Riemann surfaces, using the methods of sheaf theory that were prevailing at the time. As Brieskorn writes in his CV, this theorem was the reason why he moved from Munich to Bonn and Hirzebruch. In 1963, M. Atiyah and I. Singer generalized the Hirzebruch--Riemann--Roch theorem to elliptic differential operators on a complex manifold. This covers significant theorems of differential geometry and has important applications in theoretical physics, for which they jointly received the Abel Prize  in 2004. The theorem was further extended in a functorial way to proper morphisms of quasi-projective schemes by Grothendieck in 1967. Modifications continue to this day, e.g. with the ``Quantum Riemann--Roch'' in the Gromov--Witten theory.

The deformation theory of complex manifolds developed by Kodaira and Spencer in \cite {KS1958}, in particular the infinitesimal deformation theory, and Kuranishi's theorem on the existence of a semi--universal deformation, were later developed further and are now among the most important methods in complex analysis as well as algebraic and arithmetic geometry. I would like to mention only the existence of a semi--universal deformation for compact complex spaces, proven by Hans Grauert and Adrien Douady.

Even more important for the theory of singularities and for Brieskorn's later work was the theorem by Grauert on the existence of a semi-universal deformation of isolated singularities of complex spaces. Brieskorn himself was involved in the development of the proof by Grauert.
And that came as follows.

Brieskorn had started his professorship in G\"ottingen in July 1969, but had been on leave from September 1969 to February 1970 for a research stay at the IHES in Bures-sur-Yvette. From there he brought an interesting problem for the research seminar in G\"ottingen, jointly led by Brieskorn and Grauert in the summer and winter semester 1970/71: to prove the existence of a semi--universal deformation of isolated singularities of complex spaces. This was to be based on the papers of M. Schlessinger {\em Functors of Artin Rings} and of G. N. Tyurina {\em Locally semi - universal flat deformations of isolated singularities of complex spaces}. Tyurina's work was published in 1969 in Russian, but only in 1971 in English translation. The author of these lines, who himself attended the seminar, believes that Brieskorn himself had brought along both these papers. They were not yet known in G\"ottingen and, in particular, he had procured an English translation of Tyurina's work.

Schlessinger had specified in his work conditions for the existence of a formal semi--universal deformation, the ``Schlessinger conditions'', while Tyurina had proved the existence for normal isolated singularities, with the additional condition that the second Ext group of the holomorphic $ 1 $--forms on the singularity vanishes. Tyurina did not seem to have known Schlessinger's work, and Brieskorn's idea was that the combination of Schlessinger's and Tyurina's approaches should provide evidence of a proof for any isolated singularity without Tyurina's additional condition. He had even more precise ideas on how both works should be brought together to a proof, and he distributed the lectures accordingly to the participants of the seminary. The last lecture was assigned to Grauert with Brieskorn's comment  `and you then prove the general proposition without any assumption'. Grauert answered only  `but I need the Christmas holidays', which caused general joy.

The participation in the seminar was quite demanding, especially for someone who had just started working on his diploma thesis, but at the same time enormously stimulating and enriching. All participants were aware that they were involved in the emergence of a significant result and waited anxiously for Grauert's lecture, which was to take place in early 1971. Grauert then began his lecture with the remark that he unfortunately could not present the proof and he was not sure whether one might need an assumption like that of Tyurina. However, he could recount an interesting generalization of the Weierstrass theorem, the division with remainder by an ideal (a result found independently by Hironaka in connection with the resolution of singularities). This general division theorem then provided the crucial tool for demonstrating the existence of a semi-universal deformation of isolated singularities without any assumption that Grauert completed in the summer of 1971 (published in \cite {HG1972}).

The episode shows Brieskorn's infallible grasp of interesting and important mathematical problems, which can be seen throughout his work and in the selection of topics for diploma and doctoral theses.\bigskip

\textbf{Quotient singularities and simultaneous resolution}

The years after graduation were among the most fertile in Brieskorn's scientific life. At some point in 1963, Hirzebruch suggested that Brieskorn study and generalize the work {\em On analytic surfaces with double points} by Michael Atiyah \cite {MA1958}. In this work, Atiyah had shown, among other things, that a family $ f \! \!: \! \! X \to S $ of compact complex surfaces over a smooth $ 1 $--dimensional manifold $ S $, whose general fiber is smooth and whose finitely many special fibers only have singularities of the type $ A_1 $, has a simultaneous resolution. Here, a simultaneous resolution of a general holomorphic mapping $ f \! \!: \! \! X \to S $ is a commutative diagram of holomorphic maps
\[
\xymatrix{
Y\ar[r]^\psi\ar[d]_g & X\ar[d]^f\\
T\ar[r]^\varphi & S\ ,}
\]
where $ \varphi $ is a branched covering and $ g $ is a non--singular, proper surjective mapping that
induces for all fibers $ X_s \! \! = \! \! f ^ {- 1}  ( s) $ of $ f $ a resolution $ \psi | Y_t \! \! : Y_t \to X_s , \ \varphi (t) = s $.

Because of the local monodromy around the singular fibers of $ f $, the base change $ \varphi $ is necessary, i.e. a simultaneous resolution of $ f $ over $ S $ itself is in general not possible. The local ``geometric monodromy'' about a singular fiber $ X_ {s_0} $ can be seen as follows: restrict $ f $ to a small closed path $ \gamma $ with start and end point $ s $ around $ s_0 \in S $, so that $ \gamma $ does not go through a singular value of $ f $ and simply circles around $ s_0 $, you get a locally trivial fiber bundle over $ \gamma $, that is in general not trivial. By means of a path lifting (for example by using an Ehresmann connection) one obtains a nontrivial diffeomorphism of the fiber $ X_s $, the ``geometric monodromy'', which induces a nontrivial isomorphism of the homology of $ X_s $. The base change $ \varphi $ eliminates the local monodromy: since all the fibers of $ g $ are non-singular, the monodromy of $ g $ over $ T $ is trivial.

Hirzebruch had proposed to Brieskorn to generalize the work of Atiyah to families of surfaces with singularities of the types $ A_k, D_k, E_6, E_7, E_8 $. It turned out to be a wonderful idea, and it started a highly interesting story with many actors and great discoveries. In the end, not only the simultaneous resolution of the families of surfaces with $ ADE $ singularities was achieved, but also the discovery of exotic spheres as neighbourhood boundaries of singularities. Brieskorn writes in \cite {EB2010}
``I shall be grateful for it to my teacher until the day that I die''.

Since the problem is local, it suffices to investigate a map $ f $ from $ X = \C ^ 3 $ to $ S = \C $ (or small neighborhoods of zeros) of the form $ s = f (x, y, z) $. Here $ f (x, y, z) = 0 $ is the equation of an $ ADE $ singularity, i.e. $ f $ is a polynomial of the following list:
\[
\begin{array}{lll}
A_k: & x^{k+1}+y^2+z^2, & k\geq 1\\
D_k: & x^{k-1}+xy^2+z^2, & k\geq 4\\
E_6: & x^4+y^3+z^2\\
E_7: & x^3y+y^3+z^2,\\
E_8: & x^5+y^3+z^2\ .
\end{array}
\]

These polynomials already appeared in the works of Hermann Amandus Schwarz and Felix Klein in the 19th century (see \cite {FK1884}). Since Klein's time, they have appeared in ever new, different contexts (see \cite {GG1992} for an overview) and have fascinated mathematicians to this day. Depending on the context in which they appear, they are also called ``simple surface singularities'', ``rational double points'', ``Du--Val singularities'' or ``Kleinian singularities''.

The context in which the $ ADE $ singularities appear in Klein's work is particularly interesting to us. Klein classified the finite subgroups $ G $ from $ \SL (2, \C) $ up to conjugation and obtained the following groups:
\[
\begin{array}{ll}
C_{k+1}: & \text{the cyclic group of order } k+1\\
D_{k-2}: & \text{the binary dihedral group of order $4(k-2)$}\\
T: & \text{the binary tetrahedral group of order } $24$\\
O: & \text{the binary octahedral group of the order } $48$\\
I: & \text{the binary icosahedral group of the order } $120$
\end{array}
\]

These groups are (complex) ``reflection groups'', i.e. groups generated by complex reflections (finite--order automorphisms that fix a hyperplane), and Klein proved that the ring $ \C [z_1, z_2] ^ G $  of $ G $--invariant polynomials in $ \C [z_1, z_2] $ is generated by three invariant polynomials $ X, Y, Z $, which satisfy exactly one relation $ f (X, Y, Z) = 0 $. Klein determined the relations and found that these are given for the groups $ C_ {k + 1}, D_ {k-2}, T, O, I $ by the polynomials $ A_k, D_k, E_6, E_7, E_8 $.

Because of this result, the $ADE$ singularities are called ``2--dimensional quotient singularities''. More generally, let $ G \subset \GL (2, \C) $ be a finite subgroup acting on $ \C ^ 2 $ by matrix multiplication from the right and on $ \C [z_1, z_2] $ by $ (gf) (z_1, z_2) = f ((z_1, z_2) g)) $ for $ f \in \C [z_1, z_2] $ and $ g \in G $. The invariant ring $ \C [z_1, z_2] ^ G $ is a finitely generated $ \C $--algebra, i.e. there are finitely many invariant polynomials $ X_1, \ldots, X_n \in \C [z_1, z_2] ^ G $ with $ X_i (0) = 0 $, and finitely many relations $ f_j (X_1, \ldots, X_n) = 0 $ with $ f_j \in \C [x_1, \ldots, x_n], j = 1, \ldots, k $, so that the canonical map
\[
\C[x_1, \ldots, x_n]/\langle f_1, \ldots, f_k\rangle\to\C[z_1, z_2]^G, x_i\mapsto X_i,
\]
is an isomorphism. The bijection
\[
\C^2/G\to X:=\{x\in \C^n|f_1(x)=\cdots=f_k(x)=0\},
\]
makes the orbit space of $ G $ in a canonical way a normal analytic set in $ \C ^ n $. The space germ $ (\C ^ 2 / G, 0) = (X, 0) $ is called ``quotient singularity''. Naturally, Klein did not yet have this interpretation of the $ ADE $ singularity as quotient singularity. It is mainly due to Du Val (\cite {DV1957}).

The fact that $ADE$ singularities are quotient singularities and, of course, the explicit equations were essential for Brieskorn's proof of the simultaneous resolution of these singularities. The accomplishment of this proof was not straightforward, but was interrupted by other great discoveries of Brieskorn. In particular, the case of the $ E_8 $ singularity has caused greater difficulties, which is also reflected in the fact that Brieskorn published the proof for the $ A_k, D_k, E_6 $ and $ E_7 $ singularities in 1966 in \cite {EB1966a} but for $ E_8 $ only in 1968 in \cite {EB1968a}.

To determine the base change for a simultaneous resolution, one has to analyze the local monodromy. Since the $ADE$ singularities have weighted homogeneous (or quasihomogeneous) equations, the geometric monodromy is analytically computable, and it is of finite order. In the case of the quotient singularities of the type $ A_k, D_k, E_k $, the monodromy operation on the middle homology of the fiber is a Coxeter element of the reflection group, i.e. the product of the generators in a chamber of $ G $. So it makes sense to consider a base change $ s = t ^ d $, where $ d $ is the order of the Coxeter element, the Coxeter number. The fiber product of $ X \to S $ and the base change $ T \to S $ then has the equation
\[
f(x,y,z)-t^d=0,
\]
where $ f (x, y, z) = 0 $ is the equation of an $ADE$ singularity. In the $ A_1 $ case of Atiyah, we have the equation
\[
x^2+y^2+z^2-t^2=0,
\]
which after coordinate change has the form
\[
z_1z_2-z_3z_4=0.
\]
This is the equation of a 3-dimensional singular quadric $ Q_3 $ in $ \C ^ 4 $, i.e. the cone over a non--singular quadric in $ \P ^ 3 $. By blowing up the vertex, one obtains a non--singular variety $ Y $ over $ T $, whose exceptional divisor is isomorphic to $ \P ^ 1 \times \P ^ 1 $, which can be blown down in two ways to $ \P ^ 1 $. The resulting varieties $ Y_1 $ and $ Y_2 $ are two different simultaneous resolutions of the given family. They are so--called ``small resolutions'' of the 3-dimensional quadric $ Q_3 $; small, since the exceptional set is a (rational) curve. The natural bijective equivalence between $ Y_1 $ and $ Y_2 $ is called (Atiyah--) flop. Flops and flips play a fundamental role in the so--called minimal models of an algebraic variety, flips in the construction itself, which is known only up to dimension 3, while various minimal models are connected by a sequence of flops.

Since for $ A_k, E_6 $, and $ E_8 $ the equation $ f (x, y, z) -t ^ d = 0 $ is the form
\[
x^a + y^b + z^c + t^d = 0,
\]
Brieskorn tried to construct small modifications for these singularities by attempting to map them to other varieties for which a small resolution was already known. For the quadratic cone $ Q_3 $, this meant writing the equation $ f (x, y, z) -t ^ d $ in the form $ \phi_1 \phi_2- \phi_3 \phi_4 $. Using such methods, Brieskorn succeeded in constructing simultaneous resolutions of the $ A_k, D_k, E_6 $, and $ E_7 $ singularities in 1964 (published in \cite {EB1966a}). It also turned out that a simultaneous resolution of the map from a 3--manifold to a 1--manifold is only possible for the $ADE$ singularities. This left only the case of the $ E_8 $ singularity.

However, the $ E_8 $ singularity proved to be extremely stubborn and Brieskorn failed to construct a simultaneous resolution for them. The various attempts to do so led him to surprising results about the topology and differential topology of singularities, which I will discuss in the next section.

The simultaneous resolution of the ``icosahedron singularity'' $ E_8 $ was found by Brieskorn in September 1966 (published in \cite {EB1968a}) using very classical algebraic geometry, as he writes himself. For example, he used a paper by Max Noether from 1889 on rational dual planes and properties of exceptional curves on rational surfaces, which arise from the blowing up 8 points on a plane cubic. Brieskorn found that there are about 700 million simultaneous resolutions of $ E_8 $, exactly $ 2 ^ {14} \cdot 3 ^ 5 \cdot 5 ^ 2 \cdot 7 $, the order of the Weyl group of type $ E_8 $. The divisor class group of the  local ring of the singularity
\[
x^2+y^3+z^5+t^{30}=0
\]
has the structure of the lattice of weights of the root system of $ E_8 $. Brieskorn constructed the small resolutions of this singularity using curves with $ E_8 $ as a dual graph and thus the simultaneous resolutions of the surface singularity $ E_8 $. The various simultaneous resolutions correspond to the Weyl chambers, with the blow up of some ideal class in each chamber resulting in a simultaneous resolution.

Investigations of the simultaneous resolution of the $ADE$ singularities as quotient singularities of $ \C ^ 2 $ by a finite subgroup of $ \SL (2, \C) $ led Brieskorn to examine general quotient singularities $ \C ^ 2 / G $ in \cite {EB1968b}, where $ G $ is any finite subgroup of $ \GL (2, \C) $. He classified these singularities using results from Mumford, Hirzebruch and above all from Prill by listing all small subgroups $ G \subset \GL (2 , \C) $ (i.e. no element of $ G $ has $ 1 $ as eigenvalue with multiplicity $ 1 $). He determined the resolution graph of $ \C ^ 2 / G $, weighted by the intersection multiplicities of the exceptional curves. Brieskorn showed that this resolution graph determines the singularity up to analytic isomorphism, and from this fact he deduced the remarkable result about the uniqueness of the $ 2 $--dimensional icosahedron singularity:

\textit {The ring $ \C \{x, y, z \} / \langle x ^ 2 + y ^ 3 + z ^ 5 \rangle $ (and its completion) is the only non--regular factorial $ 2 $--dimensional analytic local ring.}

In dimension 3 there are infinitely many factorial as well as non--factorial local rings of isolated hypersurface singularities (from dimension 4 on one has always factoriality).

Brieskorn's work on simultaneous resolution and on quotient singularities played an important role in the further development of the deformation theory of rational surface singularities. I mention here only Oswald Riemenschneider, Jonathan Wahl and in connection with the program of ``minimal models'' of  Shigefumi Mori, J\'anos Koll\'ar, Miles Reid and Vyacheslav Shokurov.
\bigskip

\textbf{Topology of singularities and exotic spheres}

In September 1965 Brieskorn took up a C.L.E. Moore Instructorship at MIT in Cambridge/Massachusetts. The problem of the simultaneous resolution of $ E_8 $ had not been solved at that time and Brieskorn was looking for solutions in discussions with Heisuke Hironaka at the 1965 Arbeitstagung in Bonn and with Michael Artin and David Mumford at MIT. Brieskorn tried to compute the divisor class group of
\[
x^2+y^3+z^5+t^{30}=0,
\]
but Mumford suggested to examine first the simpler equation $ x ^ 2 + y ^ 3 + z ^ 5 + t ^ 2 = 0 $, that is, the equation of the $ 3 $--dimensional $ E_8 $ singularity. For practice, Brieskorn started with the $ 3 $--dimensional $ A_2 $ singularity
\[
z_0^3+z_1^2+z_2^2+z_3^2=0
\]
and found that it is factorial. He then re-examined the $ 3 $--dimensional $ E_8 $ singularity and, through a rather tedious explicit resolution, he showed that it is also factorial and that the second cohomology group of the singularity boundary vanishes. Since Brieskorn had expected a non-trivial divisor class group, he turned back to $ A_2 $ to better understand their topology.

Then, in September 1965, he made the unexpected discovery that  \textit {the boundary of the $ 3 $--dimensional $ A_2 $--singularity is homeomorphic to the $ 5 $--dimensional sphere.} The neighbourhood boundary of a hypersurface singularity in $ \C ^ {n +1} $ is the intersection with a sufficiently small real sphere in the $ \R ^ {2n + 2} = \C ^ {n + 1} $ around the singular point. For an isolated singularity, this is a compact $ (2n-1) $--dimensional real analytic manifold. The singularity itself, i.e. the set of zeros of the defining equation is, according to Milnor, topologically the cone over the neighbourhood boundary with the singular point as the vertex of the cone. It follows that the 3-dimensional $ A_2 $--singularity is topologically a manifold. This discovery came as a complete surprise, because in \cite {DM1961}, David Mumford had shown that isolated singularities of algebraic surfaces can never be topologically trivial, unless the singularity is analytically non--singular. Brieskorn then published in \cite {EB1966b}, that \textit {all odd $ k \geq $ 3 the singularities}
\[
z_0^3+z_1^2+\cdots +z_k^2=0,
\]
\textit{are topological manifolds}, so Mumford's theorem is a special phenomenon in dimension two.

The developments in 1965/66, which then led to the discovery of the exotic spheres as neighbourhood boundaries of singularities, are still fascinating in retrospect, above all because of the interaction of the ideas of several participants, which came about through happy circumstances. Hirzebruch reported this discovery in the seminar Bourbaki \cite {FH1967} and later at the 1996 singularity conference in Oberwolfach on the occasion of Brieskorn's 60th birthday; a short version can be found in \cite {FH1996}.

Hirzebruch reported on a conference in Rome on Brieskorn's simultaneous resolution of the singularities of types $ A_k, D_k, E_6 $, and $ E_7 $ when he received a letter from Brieskorn there on 28.09.1965, in which he wrote:

\textit {``I have made the somewhat confusing discovery in recent days that there may be $ 3 $--dimensional normal singularities that are topologically trivial. I discussed this example with Mumford this afternoon, and he has not found a mistake until this evening; here it is: $ X = \{x \in \C ^ 4 | x_1 ^ 2 + x_2 ^ 2 + x_3 ^ 2 + x_4 ^ 3 = 0 \}$}.''

This result of Brieskorn was quite exciting at that time and stimulated Hirzebruch, Milnor and others to further study the topology of isolated singularities. Of course, there was no e-mail at this time, but an extensive correspondence between Brieskorn, Hirzebruch, J\"anich, Milnor and Nash. Hirzebruch wrote to Brieskorn in March 1966 that he found a close connection between the work of Klaus J\"anich, who was also a student of Hirzebruch, on the classification of special $ O (n) $--manifolds and the neighbourhood boundary of singularities investigated by Brieskorn. J\"anich had studied the operation of a compact Lie group $ G $ on a differentiable manifold $ X $ without boundary. For special operations, the orbit space $ X / G $ is a canonically differentiable manifold with boundary. Motivated by Brieskorn's work, Hirzebruch considered the neighbourhood boundary $ \Sigma = \Sigma (k + 1, 2, \ldots, 2) $ of the $ A_k $--singularities in the $ \C ^ {n + 1} $, which was given by the following equations:
\[
\begin{array}{lcl}
z_0^{k+1}+z_1^2+\cdots +z_n^2 & = & 0\ ,\\
|z_0|^2+|z_1|^2+\cdots +|z_n|^2 & = & 1\ .
\end{array}
\]
He proved that the orthogonal group $ O (n) $ operates on $ \Sigma $  in a special way in the sense of J\"anich, with orbit space the $ 2 $-dimensional disk $ D ^ 2 $, and that  $\Sigma $ is a homology sphere for even $ k $ .

Even more exciting was his discovery that $ \Sigma $ is an exotic $ 9 $--sphere for $ n = 5 $ and $ k = 2 $, i.e. $ \Sigma $  is homeomorphic but not diffeomorphic to the standard sphere $ S ^ 9 $. The boundary $ \Sigma (3,2,2,2,2,2) $ of the $ 5 $--dimensional $ A_2 $--singularity turned out to be  the 9--dimensional exotic Kervaire--sphere, constructed by Michel Kervaire by the so-called ``plumbing'' of two copies of the tangential disk bundle of the $ 5 $--sphere.

The letter from Hirzebruch to Brieskorn dated 24.03.1966, in which he describes his discovery, was answered by Brieskorn on 29.03.1966 in the following words:

\textit { `Klaus J\"anich and I had not noticed anything about this connection of our work, and I was delighted how you brought things together.}'

While in Kervaire's construction the exotic sphere bounds a parallelizable manifold, $ \Sigma $ is the boundary of a singularity, and at first it remained mysterious where the parallelizable manifold could be found in the singularity image. Exotic spheres were first discovered by John Milnor in \cite {JM1956} and those of a fixed dimension form an abelian group $ \Theta_n $, with the connected sum as group operation. The important subgroup  $ bP_ {n + 1} $  consisting of those spheres that bound a parallelizable manifold was introduced  by Kervaire and Milnor  in 1963. They proved that the group $ \Theta_n $ is finite for $ n \geq 5 $ and that $ bP_ {4k + 2} $ is either $ 0 $ or $ \Z / 2 \Z $ and that the second case occurs exactly when the generator of $ bP_ {4k + 2} $ is the $ (4k + 1)-$dimensional exotic Kervaire sphere.

Milnor had been stimulated by Brieskorn's example of the neighbourhood boundary $ \Sigma (3,2,2,2) $ as a topological manifold to study the neighbourhood boundaries of other singularities and explained his reflections in a letter to John Nash in April 1966. Brieskorn quoted in \cite {EB2010} from this letter:

\textit{`Dear John,}

\textit{I enjoyed talking to you last week. The Brieskorn example is fascinating. After staring at it a while I think I know which manifolds of this type are spheres, but the statement is complicated and the proof doesn't exist yet. Let $\Sigma(p_1, \ldots, p_n)$ be the locus}
\[
z_1^{p_1}+\cdots + z_n^{p_m}=0\ ,\ |z_1|^2+\cdots +|z_n|^2=1\ .\text{'}
\]

Then Milnor continues with a concrete guess which of these manifolds are topological spheres. Brieskorn further mentions that the letter contains on the edge a small sketch of about 1 cm, which he would not have understood at the time.

Milnor's sketch, which I copied from \cite {EB2010} (see Fig. 1), shows the image of the Milnor fibration and thus the parallelizable manifold you are looking for. This sketch later became an icon in the theory of singularities and decorated almost every lecture on the topology of singularities.

\begin{center}
\includegraphics{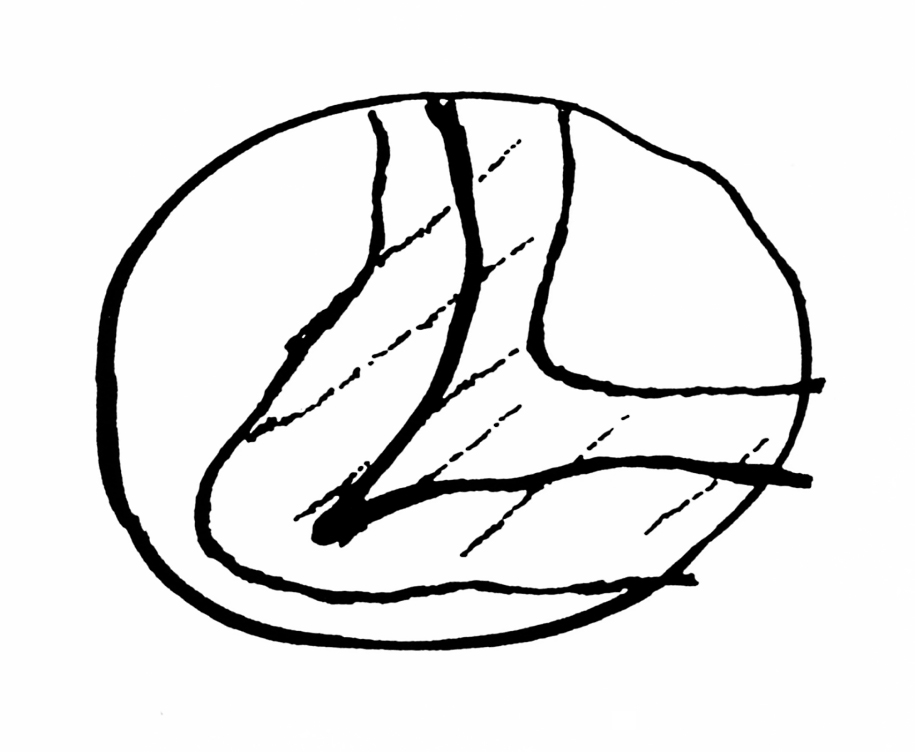}\\
Fig. 1
\end{center}

Choose a sufficiently small sphere $ B_\varepsilon $ of radius $ \varepsilon $ around the isolated singular point of the hypersurface $ f (z_1, \ldots, z_n) = 0 $ and then a small disc $ D_ \delta $ of radius $ \delta \ (<< \varepsilon) $ in the complex plane around $ 0 $. Then Milnor's sketch shows $ X = f ^ {- 1} (D_ \delta) \cap B _ {\varepsilon} $. $ X \smallsetminus f ^ {- 1} (0) $ is a differentiable fiber bundle whose non--singular fiber $ X_s = f ^ {- 1} (s), s \neq 0 $, is called ``Milnor fibre''. $ X_s $ is a $ (n-2) $--connected parallelizable manifold whose boundary is diffeomorphic to the boundary $ \Sigma $ of $ X_0 $. Thus, the Milnor fibre is the parallelizable manifold bounded by the exotic Kervaire--sphere, which Hirzebruch and Brieskorn were looking for.

Brieskorn then succeeded in \cite {EB1966c} to fully prove Milnor's conjecture within 14 days. At the same time he showed by an explicit calculation:

\textit {The neighbourhood boundary $ \Sigma (2,2,2,3,5) $ of the icosahedron singularity is Milnor's exotic $ 7 $--sphere, the creator of the group $ bP_8 = \Theta_7 $ of order 28. All the different 28 exotic differentiable structures on $ S ^ 7 $ are given by the boundary $ \Sigma (2,2,2,3,6k-1), k = 1, \ldots, 28 $, hence by simple real analytic equations. In addition, he showed that every odd dimensional sphere bounding a parallelizable manifold is diffeomorphic to the neighbourhood boundary $ \Sigma (a_1, \ldots, a_m) $.}

This was considered a sensation. While Milnor's first construction of an exotic sphere was indeed very specific, Brieskorn's construction was quite natural and anything but ``exotic''.

That Brieskorn was able to prove Milnor's conjecture so quickly is also due to a lucky circumstance. Looking through the newly published journals in MIT's library, he came across the work \cite {FP1965} by Fr\'ed\'eric Pham. Pham, motivated by the singularities of Feynman integrals in theoretical physics, examined in this paper exactly the singularities
\[
X_1^{a_1}+\cdots + X_n^{a_n}=0,
\]
which Milnor had also considered in his letter to Nash. Pham calculated for these singularities the homotopy type of the Milnor fiber and the monodromy of the Milnor fibration. Brieskorn used these results and Hirzebruch's calculation of the signature of the Milnor fibre to prove the above mentioned results. Since then, these singularities are also called ``Brieskorn Singularities'' or ``Brieskorn--Pham Singularities''.

Brieskorn's discovery of the exotic differentiable structures on the neighbourhood boundary of singularities have led to many applications in the work of other mathematicians about the differential topology of manifolds. In this context, there is only one work by Brieskorn himself, the construction of exotic Hopf manifolds, together with Antonius van de Ven in \cite {BV1968}.

Brieskorn describes the two years in Boston and Cambridge as the two best of his mathematical life.
\bigskip

\textbf{Picard--Lefschetz monodromy and Gauss--Manin connection}

For an isolated singularity, given by $f\in\C\{x_0, \ldots, x_n\}, f(0)=0$, consider Milnor's construction $f\!\!:\!\! X=f^{-1}(S)\cap B\to S$ with $B=B_\varepsilon$ and $S=D_\delta$ from the previous section. The non--singular  Milnor fibre $X_s=f^{-1}(s)$ is a deformation of $ X_0 = f ^ {- 1} (0) $, the simplest deformation since it is given by $ f $ itself. The Milnor fibre is the singularity-theoretic explanation for the fact that Brieskorn's exotic spheres bound parallelizable manifolds. However, it does not yet explain how these parallelizable manifolds can be constructed by plumbing. This requires a somewhat more complicated deformation, a so-called morsification. The idea dates back to the two volumes {\em Th\'eorie of the fonctions alg\'ebriques de deux variables ind\'ependantes} by Picard-Simart published in 1897 and 1906, and to the monograph  {\em L'analysis situs et la g\'eom\'etry alg\'ebrique} by Lefschetz. It was later developed into the local Picard--Lefschetz theory, to which Brieskorn contributed in 1970 in the appendix to \cite {EB1970a}.

There Brieskorn considers a deformation
\[
f_a(x)=f(x)-\sum\limits_{i=0}^n a_ix_i\ ,
\]
where $ a = (a_0, \ldots, a_n) $ is chosen to be sufficiently general. If $ \mu = \mu (f) $ denotes the ``Milnor number'' of $ f $, i.e. the vector space dimension of the Milnor algebra
\[
\C\{x_0, \ldots, x_n\}/\langle \frac{\partial f}{\partial x_0}, \ldots, \frac{\partial f}{\partial x_n}\rangle,
\]
then, in the neighbourhood of $ 0 \in \C ^ {n + 1} $, there are exactly $ \mu $ points $ z_r $ such that $ f_a $ has an ordinary double point  in $ z_r $, i.e. has a singularity of type $ A_1 $. A deformation with $ \mu$ ordinary double points near $ 0 $ is called a ``morsification'' of $ f $, an idea that goes back to Ren\'e Thom. The map $ f_a: X ^ a = f_a ^ {- 1} (S) \cap B \to S $ is singular exactly in the points $ z_1, \ldots, z_ \mu $  and outside the fibers through these points a differentiable  fiber bundle with fiber $ X_s ^ a = f_a ^ {- 1} (s), \ s \neq f_a (z_r), $ diffeomorphic to the Milnor fibre $ X_s $.

Brieskorn now considers Milnor's construction for the ordinary double points $ z_r $ of $ f_a $, i.e. $ f_a ^ r: f_a ^ {- 1} (D_r) \cap B_r \to D_r $, where $ B_r \subset B $ is a sufficiently small ball around $ z_r $ of radius $ \rho $ and $ D_r \subset S $ is a small disc of radius $ \delta << \rho $ around $ f_a (z_r) $. Then, for suitable coordinates, $ y_0, \ldots, y_n $ in neighbourhood of $ z_r $
\[
f_a^r(y_0, \ldots, y_n)=f_a(z_r)+y_0^2+\cdots y_n^2\ .
\]
It follows that the Milnor fibre $ F_a ^ r $ of $ f_a ^ r $ has the $ n $--dimensional sphere
\[
S_r^n=\{y|y\text{ real }, y_0^2+\cdots+y_n^2=\rho\}
\]
as a deformation retract. $ S_r ^ n \subset F_a ^ r $ defines a homology class $ d_r $ in $ H_n (F_a ^ r, \Z) $, $ r = 1, \ldots, \mu $, and these are the ``vanishing cycles'' already considered by Lefschetz, as they contract to the singular point $ z_r $ when $ \rho $ goes to $ 0 $. By choosing appropriate paths $ \gamma_r $ in $ D $ from a boundary point of $ D_r $ to the non-critical value $ s $, one can transport $ d_r $ into the Milnor fibre $ X_s ^ a $ and gets a homology class $ e_r $ in $ H_n (X_s ^ a, \Z) $. Brieskorn then shows
\[
H_n(X_s^a, \Z)=\Z e_1\oplus\cdots\oplus\Z e_\mu.
\]

This procedure provides the desired plumbing construction of the Milnor fibre for the $ADE$--singularities as follows. By choosing an Ehresmann connection
for the differentiable fiber bundle $ X ^ a \smallsetminus \underset {r} {\cup} f_a ^ {- 1} (f_a (z_r)) \to D \smallsetminus \{f_a (z_r) | r = 1, \ldots, \mu \} $
the vanishing cycles $ S_n ^ r $ themselves can be transported via $ \gamma_r $ to embedded $ n $--spheres into the Milnor fiber $ X_s ^ a $, which are also called vanishing cycles there. These vanishing cycles have tubular neighbourhoods in the Milnor fiber that are isomorphic to their tangent disc bundle. For the ADE--singularities, the vanishing cycles can be chosen in such a way that the Milnor fiber can be realized directly with the plumbing construction of these disc bundles as a parallelizable manifold.

However, Brieskorn's main goal in the paper \cite {EB1970a} was not to construct vanishing cycles by means of a morsification, but to compute the algebraic monodromy of an isolated hypersurface singularity $ f \in \C \{x_0, \ldots, x_n \} $ with $ f \! \!: \! \! X \to S $ as at the beginning of this section. The geometric monodromy is a diffeomorphism of the Milnor fiber $ X_s $ to itself given by lifting a single closed path $ \gamma $ around $ 0 $ in $ S $ with start and end point $ s $ to the total space of the fiber bundle
\[
X':=X\smallsetminus X_0\to S\smallsetminus\{0\}=:\!S'\ .
\]
The geometric monodromy induces the integral monodromy $ H ^ n (X_s, \Z) \xrightarrow {\cong} H ^ n (X_s, \Z) $ on the middle cohomology group of the Milnor fibre, the local Picard--Lefschetz monodromy of $ f $, whose characteristic polynomial $ \Delta_f $  largely determines, according to Milnor, the topology of the boundary $ \Sigma $ of $ f $.

Brieskorn gives in this paper an algebraic description of the complex local Picard--Lefschetz monodromy
\[
h_f:H^n(X_s, \C)\xrightarrow{\cong} H^n(X_s, \C)
\]
and derives from that an algorithm for computing the characteristic polynomial $ \Delta_f $.

Brieskorn uses holomorphic differential forms to compute the complex monodromy. First, the cohomology groups $ H ^ p (X_s, \C), s \in S'$, are the fibers of a holomorphic vector bundle whose sheaf of holomorphic sections is canonically isomorphic to
\[
R^nf_\ast\C_{X'}\otimes_{\C_{S'}} \ko_{S'}.
\]
Here $ R ^ nf_ \ast \C_ {X '} $ is the $ n $--th direct image sheaf of the constant sheaf $ \C_ {X'} $. Since the cohomology of the Stein manifold $ X_s $ can be calculated using holomorphic differential forms, Brieskorn considers the complex of relative holomorphic differential forms of $ X $ over $ S $,
\[
\Omega^\bullet_{X/S}=\Omega^\bullet_X/df\wedge\Omega^{\bullet-1} X,
\]
with the differential $ \Omega ^ p_ {X / S} \to \Omega ^ {p + 1} _ {X / S} $ induced by the outer derivative on the complex $ \Omega ^ \bullet_X $ of holomorphic differential forms on the manifold $ X $. One now has a canonical isomorphism
\[
R^n f_\ast\C_{X'}\otimes _{\C_{S'}}\ko_{S'} \cong H^n(f_\ast\Omega^\bullet_{X'/S'})
\]
and the right-hand side, the $ n$--th cohomology sheaf of the image sheaf complex $ f_ \ast \Omega ^ \bullet_ {X '/ S'} $, has with
\[
\kh^n(X/S):=H^n(f_\ast\Omega^\bullet_{X/S})
\]
 a continuation to all of $S$. Brieskorn shows that $ \kh ^ n (X / S) $ is coherent on $ S $ and that the stalk satisfies
 \[
\kh^n(X/S)_0=H^n(\Omega^\bullet_{X/S, 0})=:H,
\]
which depends only on the singularity of $ f $ in $ 0 $.

$ \kh ^ n (X / S) $ has as $ \ko_S $--sheaf the rank $ \mu = \mu (f) = \dim_ \C H ^ n (X_s, \C) $, $ s \in S '$, and Brieskorn conjectured that it is locally free, which was shortly thereafter proved by Marcos Sebastiani in \cite {MS1970}. Brieskorn defines on $ H $  the (meromorphic) \textit {local Gauss--Manin connection} by the formula
\[
\triangledown_f\omega=\frac{d\omega}{df}\ .
\]
This means that for a representative $ \widetilde {\omega} \in \Omega ^ n_ {X, 0} $ of $ \omega $ we have an equation $ d \widetilde {\omega} = df \wedge \psi $ where $ \frac {d \omega} {df} $ denotes the class of $ \psi $ in $ \Omega ^ n_ {X / S, 0} / d \Omega ^ {n-1} _ {X / S, 0} $. That this is well defined follows from the so--called De Rham lemma, in principle a statement about the exactness of the Koszul complex for the regular sequence $ \frac {\partial f} {\partial x_0}, \ldots, \frac {\partial f} {\partial x_n} $. For $ k $ with $ f ^ k \in \langle \frac {\partial f} {\partial x_0}, \ldots, \frac {\partial f} {\partial x_n} \rangle $ we have $ f ^ k \frac {d \omega} {df} \in H $ and Brieskorn shows:

\textit {$  \triangledown_f $ is a singular (meromorphic) first--order differential operator on $ H $ whose monodromy (by analytic continuation of a fundamental system of solutions along a closed path in $ S'$) is canonical isomorphic to the local Picard--Lefschetz  monodromy.}

In addition, Brieskorn proves that $ \triangledown_f $ is ``regular--singular'', i.e. it can be transformed by a meromorphic transformation into a differential operator with a pole of 1st order. From this, Brieskorn derives an algorithm for the computation of the characteristic polynomial $ \Delta_f $ of the monodromy of $ \triangledown_f $.

Since $ \Delta_f $ is an integer polynomial, which is algebraic in a certain sense (as Brieskorn shows), it can be deduced from the regularity of $ \triangledown_f $ together with the solution of the 7th Hilbert problem by Gelfand and Schneider (1934), that \textit {the eigenvalues of the monodromy are roots of unity $ e ^ {2 \pi i \mu_j} $ with rational $ \mu_j $}. The statement of this theorem is also referred to as ``monodromy theorem'' that had already been proved by Pierre Deligne in 1970 for global algebraic morphisms by other methods. Brieskorn's proof is considered particularly elegant.

The Manuscripta paper on the local Gauss--Manin connection led to significant developments, among others by Brieskorn's students Kyoji Saito, John Peter Scherk, Wolfgang Ebeling, Claus Hertling and the author of these lines. I myself was a student in G\"otingen and Brieskorn's first diploma student when he finished the paper. I was given the task to generalize the results to complete intersections. The main difficulty was a generalization of the lemma of De Rham. With the help of cohomological methods, which I used at the suggestion of Jean--Pierre Serre during his visit to G\"ottingen, the proof was successful and was a main result of my diploma thesis, which was completed in 1971. Later, the ``generalized De Rham--Lemma''  was further generalized by Kyoji Saito as well as Wolfgang Ebeling and Sabir Gusein--Zade. The local Gauss--Manin connection, along with Malgrange's index theorem for regular--singular differential operators, was also the key to the proof of an algebraic formula for the Milnor number of isolated complete intersection singularities in \cite {GG1975} (announced in the joint work \cite {BG1975}), which was independently derived  by L\^e D$ \tilde {\text {u}} $ng Tr\'ang using topological methods.

Besides the module $H=H^n(\Omega^\bullet_{X/S,0})$ Brieskorn introduced the two modules
\[
\begin{array}{lcl}
H': & = & df\wedge\Omega^n_{X,0}/df\wedge d\Omega^{n-1}_{X,0}\ \text{ and}\\
H'': & = & \Omega^{n+1}_{X,0}/df\wedge d\Omega^{n-1}_{X,0},
\end{array}
\]
which are also free $ \ko_{S, 0} $ modules of rank $ \mu (f) $, and the Gauss--Manin connection is then identified with a map $ \triangledown_f \! \!: \! \! H' \to H'', [df \wedge \omega] \to [d \omega] $. The meaning of this ad hoc definition was not clear at first, but it later turned out to be fundamental. $ H'$ and $ H''$ are today referred to as the ``Brieskorn--lattice'' and especially $ H''$ plays an important role in the study of the mixed Hodge structure of isolated singularities and in Kyoji Saito's ``higher--residue pairings''. In addition to the already mentioned students of Brieskorn important works in this context are due to Morihiko Saito, Daniel Barlet, Claude Sabbah, Mathias Schulze and Christian Sevenheck,  to name but a few.
\medskip

\textbf{Simple singularities and simple Lie groups}

Brieskorn's work on simultaneous resolution of simple singularities led to one of his most important findings, the relationship between $ADE$ singularities and simple Lie groups. He reported on this discovery to the International Congress of Mathematicians in Nice in 1970 and published the result in the short note \cite {EB1970b}.

Alexander Grothendieck had read Brieskorn's work on simultaneous resolution and was led to a conjecture which he told Brieskorn. While Brieskorn had studied $ 1 $--parametric deformations of the ADE singularities given by the defining polynomial, Grothendieck proposed to look at the semiuniversal deformation of these singularities. He suggested that this is determined by the adjoint quotient map of the simple Lie algebra of type $ A $, $ D $ or $ E $ and that a simultaneous resolution of the semiuniversal deformation of the singularities of the corresponding type is given with the hep of  the Springer--resolution of the nilpotent variety. Grothendieck himself had studied simultaneous resolutions of singularities of adjoint quotient maps and had come across the suspected connection through Brieskorn's work.

Let $ G $ be a simple complex (algebraic) Lie group, i.e. a complex algebraic manifold with regular group action that is simple as a group. If $ G $ is simply connected, then $ G $ is uniquely determined by its Lie algebra $ \mathfrak {g} $ up to isomorphism. The simple Lie groups correspond to the simple Lie algebras and these are classified by their fundamental root system. The root systems, in turn, are described by their ``Dynkin--diagram'' (also Coxeter--Dynkin--Witt diagram) and determine $ G $ up to isomorphism.
The classification of all Dynkin--diagrams resulting from the simple Lie groups provides four infinite series $A_k(k\geq 1), B_k (k\geq 2), C_k (k\geq 3), D_k (g\geq 4)$ and the five exceptional cases $E_6, E_7, E_8, F_4$ and $G_2$. The Dynkin--diagrams of type $ A_k, D_k, E_6, E_7, E_8 $ are characterized by the fact that they are homogeneous, i.e. their root systems have equal roots. These diagrams have the following shape ($ADE$ graphs):

%

$A_k: \includegraphics{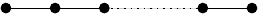}$\hskip1cm $D_k:\includegraphics{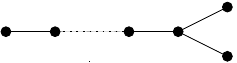}$\ \ \ \hskip 1cm ($k$ vertices)\\

$E_6: \includegraphics{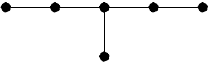}$\hskip1.55cm$E_7: \includegraphics{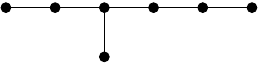}$\hskip 1cm  $E_8: \includegraphics{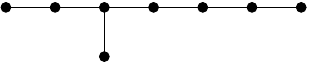}$

The name $ADE$ singularity for the quotient singularities of the finite subgroups of $ \SL (2, \C) $ comes from the relation to the simple Lie groups of type $ A_k, D_k $ or $ E_k $. At the ICM 1970 in Nice, Brieskorn presented the construction of the $ADE$ singularities and their seminuniversal deformation directly with the help of the corresponding Lie group as follows.

Considering the operation of the simple complex Lie group $ G $ on itself by conjugation, we call $ x \in G $ ``regular'' if the orbit of $ x $, that is its conjugate class in $ G $, has the maximum dimension. If $ d $ is this dimension, then the next smaller orbit dimension is $ d-2 $ and elements of this orbit dimension are called ``subregular''. $ G $ has exactly one regular orbit, containing $ 1 \in G $ in its closure, and the closure of this orbit is called Uni$ (G) $, since it is the variety of the unipotent elements of the group. The complement of the regular orbit in Uni$ (G) $ has codimension 2 and is itself the closure of exactly one subregular orbit. If $ x \in \text {Uni} (G) $ is an arbitrary element, then consider a small slice $ X \subset G $ through $ x $ transversal to the orbit of $ X $ and a regular projection $ \pi: (G, x) \to (X, x) $ of complex space germs. The space germ $ (X, x) $ has in $ (G, x) $ complementary dimension to the orbit of $ x $ and is smooth if $ x $ is a regular element. Only slices through non-regular orbits produce singularities.
If $ x $ is subregular, $ X \cap \text {Uni}(G) $ has dimension two and an isolated singularity in $ x $.

Let $ x = x_sx_n $ be the Jordan--decomposition of $ x \in G $ into a semisimple and unipotent part. Assigning to $ x $ the conjugate class of $ x_s $ yields a morphism $ \Phi: G \to T / W $, where $ T $ is a maximal torus in $ G $ and $ W $ is the Weyl group. $ \Phi $ is the adjoint quotient map. Each fiber of $ \Phi $ is the union of finitely many conjugation classes, $ T / W $ is a $ k $--dimensional complex manifold ($ k $ = number of vertices of the Dynkin diagram) and $ \Phi $ maps Uni$(G) $ to $ 1 \in T / W $. With these notations Brieskorn proved the following in \cite {EB1970b}.

Let $ G $ be a simply connected complex Lie group of the type $ A_k, D_k, E_6, E_7, E_8 $ and $ x \in G $ a subregular unipotent element. Then:

\begin{enumerate}
\item [(1)] \textit{$(X\cap \text{ Uni }(G),x)$ is isomorphic to an ADE singularity of the same type as $G$.}

\item [(2)] \textit{The adjoint quotient map germ in $x$ factorises as $\Psi\circ\pi$,
\[
\Phi: (G,x)\xrightarrow{\pi}(X,x)\xrightarrow{\Psi}(T/W, 1),
\]
where $\Psi$ is the semi--universal deformation of the corresponding quotient singularity.}
\item [(3)] \textit{Let
\[
\xymatrix{\Gamma\ar[r]\ar[d] & G\ar[d]\\
T\ar[r] & T/W}
\]
be the simultaneous resolution of the adjoint quotient map of Grothendieck, with $\Gamma=\{(x,B)|x\in G, B \text{ Borel--subgroup containing } x\}$ and $Y$ the preimage of the transversal slice  $X$ in $\Gamma$. Then
\[
\xymatrix{Y\ar[r]\ar[d] & X\ar[d]\\
T\ar[r] & T/W
}
\]
is a simultaneous resolution of the semiuniversal deformation of the quotient singularity of type  $A_k, D_k, E_k$.}
\end{enumerate}

Brieskorn's proof  makes essential use of the fact that $ \Phi $ is given by weighted-homogeneous polynomials and that the $ADE$ singularities are characterized by their weights. Incidentally, the Weyl group $ W $ is equal to the monodromy group of the singularity. The proof sketched by Brieskorn in \cite {EB1970b} was completely worked out by Peter Slodowy. Slodowy later extended the construction to all simple Lie groups, including the inhomogeneous root systems $ B_k, C_k, F_4 $ and $ G_2 $, even over fields of arbitrary characteristic \cite {PS1980}.

The relationship between Lie groups and singularities shown by Brieskorn led to further investigations. An entirely different construction of the $ADE$ singularities with the help of the simple algebraic groups of type $ A, D $ and $ E $ is due to Friedrich Knop \cite {FK1987}, although the singularities there are realized in different dimensions. A clarification of the occurrence of the polyhedral groups in Brieskorn's construction, and thus a direct relationship between the simple Lie groups and the finite Klein groups, was achieved by Peter Kronheimer \cite {PK1990}, using differential geometric methods. Brieskorn had still written at the end of \cite {EB1970b}:

{\em `Thus we see that there is a relation between exotic spheres, the icosahedron and $E_8$. But I still do not see why the regular polyhedra come in.}'

Slodowy also investigated more complicated singularities and associated them with Kac--Moody Lie algebras. Shortly before his death in 2002, he succeeded in constructing all the simply elliptic singularities with the aid of the adjoint quotient map of the infinite--dimensional loop group. The work was completed by Stefan Helmke in \cite {HS2004}. Slodowy had searched for this result for many years and he reported this to Brieskorn, a few days before his death when he was already badly marked by his illness, but still full of passion and enthusiasm. Brieskorn was deeply touched that the continuation of his ideas could bring consolation and joy even in the hardest hour.
\bigskip

\textbf{Generalized braid groups, Milnor lattice and Lorentzian space--forms}

The construction of the semi-universal deformation of an ADE singularity with the help of the adjoint operation of the corresponding simple Lie group led Brieskorn to investigate operations of generalized braid groups and thus to turn to the investigation of discrete structures of isolated singularities.

Let $ W $ be a finite reflection group operating linearly on the real finite--dimensional Euclidean vector space $ E'$ and let $ D' \subset E'$ be the union of the reflection hyperplanes $ H'_s $, where $ s $ is an element of the set of the reflections $ \Sigma $ in $ W $. Brieskorn considers the complexification  $ E $ of $ E'$ resp. $ H_s $ of $ H'_s $ and the union $ D \subset E $ of the $ H_s, s \in \Sigma $. The operation of $ W $ on $ E'$  extends canonically to $ E $ and maps $ E_ {reg} = E \smallsetminus D $ to itself. $ E_ {reg} / W $ is the space of regular orbits of the finite complex reflection group $ W $, whose fundamental group was calculated by Brieskorn in  \cite {EB1971} (see also \ cite {EB1973}). He shows:

\textit{The fundamental group $\Pi_1(E_{reg}/ W)$ has a presentation with generators $g_s, s\in\Sigma$, and relations}
\[
g_sg_tg_s\cdots=g_tg_sg_t\cdots ,
\]
\textit{with $m_{st}$ factors on both sides.}

Here $(m_{st})$ is the Coxeter--matrix of $W$ with $m_{st}=$ order of $st$, and the $ g_s $ are given by an explicit geometric construction.

For the symmetric group $ W \! \! = \! \! S_n \ (= A_ {n-1}) $, the corresponding fundamental group is the braid group $ B_n$ introduced in 1925 by Emil Artin, the father of Michael Artin, as proved by Fox and Neuwirth in 1962 and independently by Arnold in 1968. The finite irreducible reflection groups are classified and fall into the types $ A_k, B_k, D_k, E_6, E_7, E_8, F_4, G_2, H_3, H_4 $ and $ I_2 (m), m = 5 $ or $ m \geq 7 $. The fundamental groups of the regular orbits of these complex reflection groups are therefore generalizations of the braid groups.

The connection with singularities comes from the fact that for $ W $ of type $ A_k$, $D_k$, $E_6$, $E_7$, $E_8 $ the space $ E_ {reg} / W $ is the complement of the discriminant in the base space of the semiunversal deformation of the simple singularity of the same type. This follows from Brieskorn's construction in \cite {EB1970b}.

These generalized braid groups were baptized ``Artin groups'' by Brieskorn and Kyoji Saito in honor of Artin in \cite {BS1972}, and they examined them from a combinatorial point of view. Among other things, they solve the word and conjugation problems for these groups and determine the center. These results were obtained at about the same time by Pierre Deligne in \cite {PD1972}.
Deligne proved that the spaces $ E_ {reg} / W $  considered above are Eilenberg--MacLane spaces, as was conjectured by Brieskorn in \cite {EB1971}.

In the works mentioned below, Brieskorn studies discrete invariants of special classes of singularities. Let $ (X_0, x) \subset (\C ^ {n + 1}, x) $ be an isolated hypersurface singularity and $ F \! \!: \! \! X \to S, F (x) = 0 $, a suitable representative of the semiuniversal deformation of $ (X_0, x) $. If $ D $ is the discriminant of $ F $, that is, the set of points $ s \in S $ for which the fiber is not smooth, then, with $ S'= S \smallsetminus D $, the restriction $ F \!: \! X'= F ^ {- 1} (S') \to S '$ is a differentiable fiber bundle with fiber $ X_s $, diffeomorphic to the Milnor fiber of $ (X_0, x) $. Since $ X_s $ has the homotopy type of a bouquet of $ n $--dimensional spheres, the middle homology group $ H_n (X_s, \Z)$ is free of rank $ \mu $, the Milnor number of $ (X_0, x) $. If $ n $ is even, $ H_n (X_s, \Z) $ carries an integral symmetric quadratic form, the intersection form $ <,> $, and the integer lattice
\[
L=H_n(X_s, \Z)
\]
is called the ``Milnor lattice'' of the singularity.

If one selects a generic complex line near $ 0 $ in the affine space containing  $ S $, then the intersection with $ S $ is a small disc $ \Delta $ intersecting the discriminant $ D $ in $ \mu $ different points $ c_1 , \ldots, c_\mu $. The restriction of $ F $ over $ \Delta $ is a morsification of $ (X_0, x) $, as described above. For $ s \in \Delta'= \Delta \smallsetminus \{c_1, \ldots, c_\mu \} $ and a choice of paths $ \gamma_i $ in $ \Delta' $ from $ s \in \Delta'$ to points near the $ c_i, i = 1, \ldots, \mu $, one obtains so--called ``vanishing cycles'' $ e_i \in H_n (X_s, \Z) $ with $ \langle e_i, e_i \rangle = -2 $. The set of vanishing cycles is denoted $ \Delta ^ \ast \subset L $.

By a proper choice of the paths $ \gamma_i $ the $ e_1, \ldots, e_\mu $ form a basis of the lattice $ L $, which is then called a ``distinguished basis''. The set of all distinguished bases of $ L $ is denoted $ B ^ \ast $.

For every basis $ B \in B ^ \ast $, the matrix of scalar products of the basis elements describes the bilinear form on $ L $, which is characterized by a graph $ D_B $. The vertices $ \{1, \ldots, \mu \} $ of $ D_B $ correspond to the basis elements $ e_1, \ldots, e_\mu $ and two vertices $ i, j $ are connected by $ | \langle e_i, e_j \rangle | $ edges, each with the sign $ \pm 1 $ of $ \langle e_i, e_j \rangle \in \Z $. $ D_B $ is called (Coxeter--) Dynkin--diagram of $ B $ and the set of all Dynkin--diagrams is denoted by $ D ^ \ast $.

On $ B ^ \ast $ and thus on $ D ^ \ast $ exists a natural operation of the classical braid group $ B_\mu $ with $ \mu $ strands, which can be described by elementary operations at the level of the paths. Brieskorn points out at various places that understanding this operation should be essential for an understanding of the semiuniversal deformation.

Another invariant is the ``monodromy group'' of $ (X_0, x) $. By definition, this is the image under the canonical homomorphism of the fundamental group of the complement of the discriminant in the automorphism group of the lattice $ L $. It is already generated by the automorphisms that belong to a distinguished basis.

An overview of these invariants and the relationships between them is given by Brieskorn in the survey article \cite {EB1983a}, and he stresses their importance for the understanding of the geometry of the semiuniversal deformation. Important work on these invariants are the fundamental works of Andrei Gabrielov \cite {AG1974} and Sabir Gusein--Zade \cite {SG1977} as well as the lecture notes of Wolfgang Ebeling \cite {WE1987}.

A first step is to understand the deformation relations between singularities of a fixed modality class, for, if one singularity deforms into another, this induces an inclusion of the corresponding Milnor lattices. The classification of isolated hypersurface singularities in terms of their modality (i.e., the number of independent parameters (moduli) of isomorphism classes in a neighbourhood of the origin of the semiuniversal deformation) was initiated by V.I. Arnold in \cite {VA1972} and is one of the starting points of singularity theory with far--reaching results. The adjacencies (deformation relations) between the $ADE$ singularities were already determined by Arnold. In \cite {EB1979} Brieskorn calculates all possible adjacencies within Arnold's list of unimodular singularities, which is the next more complicated class in Arnold's hierarchy, after the $ADE$ singularities. This work is refined in \cite {EB1983b}, where Brieskorn gives a very detailed description of the Milnor lattice of the 14 exceptional unimodular singularities.

The deformation relations within the unimodular singularities have been linked by Brieskorn with a theory that seems to be far away from the theory of singularities, namely the theory of partial compactifications of bounded symmetric domains. If $ \kf (L) $ denotes the isotropy--flag--complex of $ L $, a building in the sense of Tits, then the monodromy group $ \Gamma $ operates on $ \kf (L) $ and the $ 1 $--dimensional simplicial complex $ \kf (L) / \Gamma $ is finite. For the simplest exceptional singularities $ E_{12}, Z_{11}, Q_{10} $, Brieskorn proves that the Baily--Borel compactification of $ \kf (L) / \Gamma $ can be identified with the $ \C ^\ast $--quotient of the punctured negatively graded part of the base space of the semiuniversal deformation. Independent of Brieskorn, Eduard Looijenga proved these results for all triangular singularities $ T_ {p, q, r} $ in \cite {EL1983} and later, in a more general context, he constructed important new compactifications of locally symmetric varieties.

In the work \cite {EB1988} Brieskorn gives an overview of the operation of the braid group on the set $ B ^ \ast $ of distinguished bases of an isolated singularity. He also introduces the concept of an automorphic set, which unifies many aspects of the braid group operation, and which was later taken up and generalized by several authors. The following quote from the introduction shows again Brieskorn's joy in the unity of mathematics, which is expressed in the interplay of many different areas of mathematics.

\textit{``The beauty of braids is that they make ties between so many different parts of mathematics, combinatorial theory, number theory, group theory, algebra, topology, geometry and analysis, and, last but not least, singularities.}''

This brings me almost to the end of the review of Brieskorn's mathematical work.  Still to mention is the textbook {\em Plane algebraic curves} written together with Horst Kn\"orrer, whose latest reprint appeared in English translation in 2012. And of course his two textbooks {\em Lineare Algebra und Analytische Geometrie I, II} \ \footnote {A third volume {\em Linear Algebra and Analytical Geometry III} was only partially completed by Brieskorn. These parts with valuable historical information and cross--links are to be put online soon.} (Linear Algebra and Analytic Geometry I, II), which are also worth reading because of the historical remarks by Erhard Scholz.

There is also a mathematical work together with his students Anna Pratoussevitch and Frank Rothenh\"ausler from the year 2003 \cite {EB2003}, which I would like to mention. The origins of this work date back to at least 1992, when Brieskorn's student Thomas Fischer discovered a polyhedron that in a sense generalizes the classical dodecahedron. It has a very similar combinatorial structure to the dodecahedron, but with an axis of symmetry of order 7 instead of 5. Let $ \Gamma $ be a discrete subgroup of the Lie group $ \widetilde {SU} (1,1) $, which operates by left translations. The quotient $ \Gamma \backslash \widetilde {SU} (1,1) $ is a ``Lorentzian space--form'' and is described by a fundamental domain $ F $ for $ \Gamma $. For co-compact $ \Gamma $, according to results of Dolgachev, $ \Gamma \backslash \widetilde {SU} (1,1) $ is the boundary of a quasi-homogeneous surface singularity. The quasi-homogenous singularities of Arnold's series $ E_k, Z_k $ and $ Q_k $ are of this type.

In this work, the authors describe the fundamental domains for the corresponding groups $ \Gamma $ as polyhedra with total geodesic faces in the 3--dimensional Lorentz space, each series showing a regular characteristic combinatorial pattern associated with the classical polyhedra.

Brieskorn describes in the movie {\em Science Lives: Egbert Brieskorn}, see \cite {EB2010}, the great joy he felt when Fischer discovered the $ E_{12} $ polyhedron and as Pratoussevitch could expand this to the infinite series $ E_k, Z_k $ and $ Q_k $. In the sequence ``Melencolia'' of the film he explains that the correct beginning of this infinite series should be the polyhedron in D\"urer's famous etching \textit {Melencolia I}. In the same sequence of the movie, he also discusses the importance of intuition, especially in teaching students, in contrast to a purely analytical and structural approach. Fischer's discovery pleased Brieskorn so much that he himself calculated and drew a graphic representation of it and called it ``Opus 2''. He commissioned an artist to make a  brass 3D--sculpture of the $ E_ {12} $ polyhedron, which he donated to his teacher Friedrich Hirzebruch on his 75th birthday. As far as I know, he made only two or three copies of which he gave me a copy, which made me extremely happy. A picture of this specimen, which closes the circle from the beginnings of the Platonic solids and quotient singularities to the Lorentzian space--forms, may be a fitting conclusion to this review of Brieskorn's mathematical work.

\begin{center}
\includegraphics[width=6cm,height=6cm]{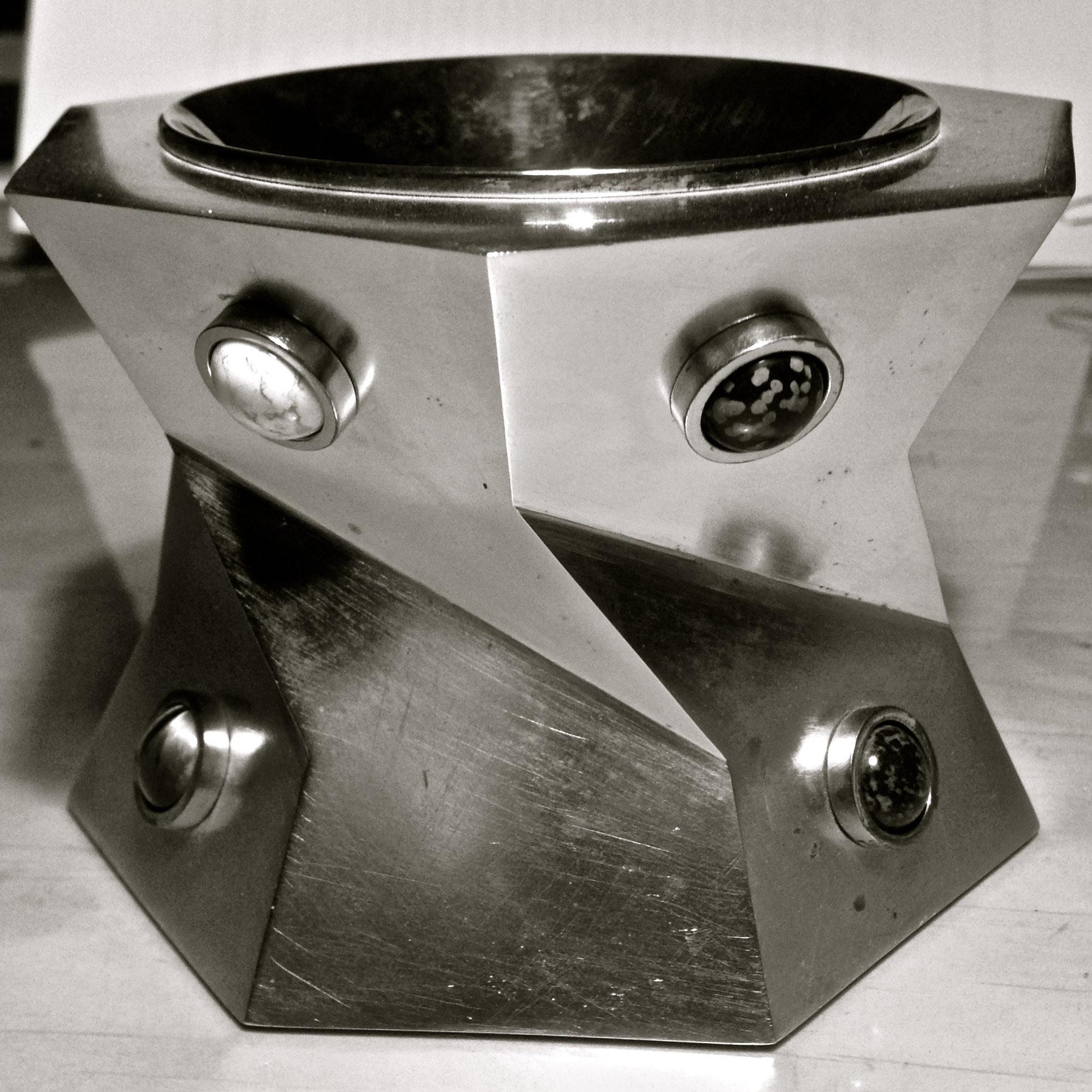}\\
$E_{12}$--Polyeder
\end{center}

Everyone who knew Egbert Brieskorn valued his extensive mathematical knowledge, his immensely broad general education, his keen intelligence, his straightforwardness and absolute intellectual honesty, his kind attention and helpfulness, and his prudent advice. And everyone who got to know him better knows that he was a kind--hearted person. The history of science will preserve his work and his name, and those who know and appreciate him will not forget him.
\bigskip

\textbf{\Large Habilitations}

Brieskorn supervised 24 Ph.D. dissertations, which can be found in the ``Mathematics Genealogy Project''. Of his doctoral students, seven have completed their habilitation:

1975:  HAMM, HELMUT AREND: Zur analytischen und algebraischen Beschreibung  der  Picard-Lefschetz-Monodromie. Göttingen.

1980:  GREUEL, GERT-MARTIN: Kohomologische Methoden in der Theorie
	isolierter Singularitäten. Bonn.

1984:  SLODOWY, PETER: Singularitäten, Kac-Moody-Lie-Algebren.
	assoziierte Gruppen und Verallgemeinerungen. Bonn.

1985:   KNÖRRER, HORST: Geometrische Aspekte integrabler Hamiltonscher
	Systeme. Bonn.

1986:  EBELING, WOLFGANG: Die Monodromiegruppen der isolierten 	
	Singularitäten vollständiger Durchschnitte.Bonn.

1986:  SCHOLZ, ERHARD: Symmetrie – Gruppe – Dualität. Studien zur
	Beziehung zwischen theoretischer Mathematik und Anwendungen
	in Kristallographie und Baustatik im 19. Jahrhundert. Wuppertal.

2000:  HERTLING, CLAUS: Frobenius-Mannigfaltigkeiten, Gauß-Manin-
	Zusammenhänge und Modulräume von Hyperflächensingularitäten.
	Bonn.

\renewcommand{\refname}{Literature}

\bigskip

G.-M. Greuel, Technische Universit\"at Kaiserslautern, Kaiserslautern, Germany\\
e-mail: greuel@mathematik.uni-kl.de
\medskip

W. Purkert, Universit\"at Bonn, Bonn, Germany\\
e-mail:  edition@math.uni-bonn.de

\end{document}